\documentclass[12pt]{amsart}

\usepackage[matrix,arrow,curve,frame]{xy}
\usepackage{amsmath,amsthm,amssymb,enumerate}
\usepackage{latexsym}
\usepackage{amscd}
\usepackage[colorlinks=false]{hyperref}
\usepackage{euscript}

\newtheorem{thm}[subsection]{Theorem}
\newtheorem{defn}[subsection]{Definition}

\newtheorem{lemma}[subsection]{Lemma}

\newtheorem{remark}[subsection]{Remark}

\theoremstyle{definition}

\def\cO{{\cal O}}

\def\cA{{\cal A}}

\def\ra{\rightarrow}
\def\bra{\langle}
\def\ket{\rangle}

\def\H{{\bf H}}

\def\L{{\bf L}}

\def\cA{{\mathcal A}}

\def\cC{{\mathcal C}}
\def\cD{{\mathcal D}}
\def\cE{{\mathcal E}}
\def\cF{{\mathcal F}}

\def\cI{{\mathcal I}}

\def\cK{{\mathcal K}}

\def\cO{{\mathcal O}}

\def\cQ{{\mathcal Q}}
\def\cR{{\mathcal R}}
\def\cS{{\mathcal S}}

\def\cV{{\mathcal V}}
\def\cW{{\mathcal W}}

\def\gg{{\mathfrak g}}

\def\gl{{\mathfrak l}}

\def\gs{{\mathfrak s}}

\newfont{\german}{eufm10}

\begin{document}
\pagestyle{plain}

\title
{Chiral equivariant cohomology of a point: a first look}

\author{Andrew R. Linshaw}
\address{Fachbereich Mathematik, Technische Universit\"at Darmstadt, 64289 Darmstadt, Germany.}
\email{linshaw@mathematik.tu-darmstadt.de}
%\thanks{}

%\date{\today}
{\abstract

\noindent
The chiral equivariant cohomology contains and generalizes the classical equivariant cohomology of a manifold $M$ with an action of a compact Lie group $G$. For any simple $G$, there exist compact manifolds with the same classical equivariant cohomology, which can be distinguished by this invariant. When $M$ is a point, this cohomology is an interesting conformal vertex algebra whose structure is still mysterious. In this paper, we scratch the surface of this object in the case $G=SU(2)$.}

\maketitle
\tableofcontents
\section{Introduction}

\noindent

Let $G$ be a compact, connected Lie group with complexified Lie algebra $\gg$, and let $M$ be a topological space on which $G$ acts by homeomorphisms. The equivariant cohomology $H^*_G(M)$ is defined to be $H^*((M\times E)/G))$, where $E$ is any contractible space on which $G$ acts freely. This is known as the {\it Borel construction}. If $M$ is a smooth manifold, and the action of $G$ on $M$ is by diffeomorphisms, there is a de Rham model of $H^*_G(M)$ due to H. Cartan, and developed further by Duflo-Kumar-Vergne \cite{DKV} and Guillemin-Sternberg \cite{GS}. In fact, one can define the equivariant cohomology $H^*_G(A)$ of any $G^*$-algebra
$A$. A $G^*$-algebra is a commutative superalgebra equipped with an action of $G$, together with a compatible action of a certain differential Lie superalgebra $(\gs\gg,d)$ associated to the Lie algebra $\gg$ of $G$. The de Rham model of $H^*_G(M)$ is obtained by taking $A$ to be the algebra $\Omega(M)$ of differential forms on $M$, and $H^*_G(\Omega(M)) \cong H^*_G(M)$ by an equivariant version of the de Rham theorem.

In \cite{LLI}, the chiral equivariant cohomology $\H^*_G(\cA)$ of an $O(\gs\gg)$-algebra $\cA$ was introduced as a vertex algebra analogue of the equivariant cohomology of $G^*$-algebras. Examples of $O(\gs\gg)$-algebras include the semi-infinite Weil complex $\cW(\gg)$, which was introduced by Feigin-Frenkel in \cite{FF}, and the chiral de Rham complex $\cQ(M)$ of a smooth $G$-manifold $M$, which was introduced by Malikov-Schechtman-Vaintrob in \cite{MSV}. In \cite{LLSI}, the chiral equivariant cohomology functor was extended to the larger categories of $\gs\gg[t]$-algebras and $\gs\gg[t]$-modules. The main example of an $\gs\gg[t]$-algebra which is not an $O(\gs\gg)$-algebra is the subalgebra $\cQ'(M)\subset \cQ(M)$ generated by the weight zero subspace. 
$\H^*_G(\cQ(M))$ and $\H^*_G(\cQ'(M))$ are both \lq\lq chiralizations" of $H^*_G(M)$, that is, vertex algebras equipped with weight gradings
$$\H^*_G(\cQ(M)) = \bigoplus_{m\geq 0} \H^*_G(\cQ(M))[m],\ \ \ \ \ \H^*_G(\cQ'(M)) = \bigoplus_{m\geq 0}\H^*_G(\cQ'(M))[m],$$ such that $\H^*_G(\cQ(M))[0] \cong H^*_G(M) \cong  \H^*_G(\cQ'(M))[0]$. If $G$ acts on $M$ effectively, $\H^*_G(\cQ(M))$ is purely classical, and $\H^*_G(\cQ(M))[m]$ vanishes for $m>0$ (see Theorem 1.5 of \cite{LLSII}). The functor $\H^*_G(\cQ'(-))$ is more interesting, and is a stronger invariant than $H^*_G(-)$ on the category of compact $G$-manifolds. For any simple $G$, there exists a sphere $S$ with infinitely many smooth actions of $G$ with the same classical equivariant cohomology, which can all be distinguished by $\H^*_G(\cQ'(S))$ (see Theorem 1.7 of \cite{LLSII}).

In the case where $M$ is a point $pt$, $\cQ(pt) = \cQ'(pt) = \mathbb{C}$. We will refer to $\H^*_G(\mathbb{C})$ as the {\it chiral point algebra}; it plays the role of $H^*_G(pt) \cong Sym(\gg^*)^G$ in the classical theory, and $\H^*_G(\mathbb{C})[0]\cong H^*_G(pt)$. For any $\gs\gg[t]$-module $\cA$, there is a chiral Chern-Weil map \begin{equation}\kappa_G: \H^*_G(\mathbb{C}) \ra \H^*_G(\cA).\end{equation} When $\cA = \cQ'(M)$ for some $G$-manifold $M$, this extends the classical Chern-Weil map at weight zero. We regard the elements of $\H^*_G(\mathbb{C})$ as universal characteristic classes, and an important problem in this theory is to describe $\H^*_G(\mathbb{C})$ for any $G$.

If $G$ is abelian, $\H^*_G(\mathbb{C})$ is the abelian vertex algebra generated by $H^*_G(pt)$, but for nonabelian $G$, $\H^*_G(\mathbb{C})$ possesses a rich algebraic structure. For semisimple $G$, $\H^*_G(\mathbb{C})$ has a conformal structure of central charge zero, and the Virasoro class $\L$ has no classical analogue in $H^*_G(pt)$. In this paper, we focus on the simplest nontrivial case $G=SU(2)$. Our main result is the construction of an injective linear map $$\Psi: U(\gs\gl_2)\rightarrow \H^*_{SU(2)}(\mathbb{C})$$ whose image consists entirely of nonclassical elements. Let $x,y,h$ be the usual root basis for $\gs\gl_2$, satisfying $[x,y]=h$, $[h,x]=2x$, and $[h,y] = -2y$. Given a monomial $\mu = x^r y^s h^t\in U(\gs\gl_2)$, $\Psi(\mu)$ is homogeneous of degree $4r-4s$ and conformal weight $2s+t+2$. In particular, $\Psi$ maps the eigenspace of $[h,-]$ of eigenvalue $d$ into $\H^{2d}_{SU(2)}(\mathbb{C})$. Moreover, the Virasoro element $\L$, which has degree zero and weight two, is precisely $\Psi(1)$. We conjecture that the image of $\Psi$, together with the classical generator of $H^*_{SU(2)}(pt)$, forms a strong generating set for $\H^*_{SU(2)}(\mathbb{C})$. (In this terminology, a vertex algebra $\cA$ is said to be strongly generated by a set $S$, if $\cA$ is spanned by the set of iterated Wick products of elements of $S$ and their derivatives). For the sake of illustration, we write down explicit representatives for a few of the classes given by $\Psi$, and we compute some relations among them. 

The main technical difficulty in studying the chiral point algebra is that the complex which computes this cohomology is a commutant subalgebra of a $bc\beta\gamma$ system, and commutant vertex algebras of this kind are not well understood. Much of this paper is devoted to developing techniques for describing such commutant algebras. Our main tool is the notion of a {\it good increasing filtration} on a vertex algebra, which was introduced in \cite{LiII} and used in \cite{LL} to study the commutant problem. The associated graded object of such a vertex algebra is always a supercommutative ring with a derivation. By passing to the associated graded object, one can apply techniques from commutative algebra and algebraic geometry, in particular the theory of jet schemes.

An interesting problem which we do not address in this paper is to give an alternative, more geometric construction of the chiral equivariant cohomology, which is suitable for any topological $G$-space $M$, and gives the same cohomology as our theory when $M$ is a manifold. In other words, we would like to find a chiral analogue of the Borel model for equivariant cohomology. Such a construction would necessarily include a topological realization of the chiral point algebra. It is our hope that the structure we describe in this paper in the case $G=SU(2)$ may give some hint about where to look for such a construction. 

\section{Vertex algebras}
In this section, we define vertex algebras, which have been discussed from various different points of view in the literature (see for example \cite{B}\cite{FLM}\cite{K}\cite{FBZ}). We will follow the formalism developed in \cite{LZI} and partly in \cite{LiI}. Let $V=V_0\oplus V_1$ be a super vector space over $\mathbb{C}$, and let $z,w$ be formal variables. By $QO(V)$, we mean the space of all linear maps $$V\rightarrow V((z)):=\{\sum_{n\in\mathbb{Z}} v(n) z^{-n-1}|
v(n)\in V,\ v(n)=0\ \text{for} \ n>>0 \}.$$ Each element $a\in QO(V)$ can be
uniquely represented as a power series
$$a=a(z):=\sum_{n\in\mathbb{Z}}a(n)z^{-n-1}\in (End\ V)[[z,z^{-1}]].$$ We
refer to $a(n)$ as the $n$th Fourier mode of $a(z)$. Each $a\in
QO(V)$ is assumed to be of the shape $a=a_0+a_1$ where $a_i:V_j\ra V_{i+j}((z))$ for $i,j\in\mathbb{Z}/2\mathbb{Z}$, and we write $|a_i| = i$.

On $QO(V)$ there is a set of nonassociative bilinear operations
$\circ_n$, indexed by $n\in\mathbb{Z}$, which we call the $n$th circle
products. For homogeneous $a,b\in QO(V)$, they are defined by
$$
a(w)\circ_n b(w)=Res_z a(z)b(w)~\iota_{|z|>|w|}(z-w)^n-
(-1)^{|a||b|}Res_z b(w)a(z)~\iota_{|w|>|z|}(z-w)^n.
$$
Here $\iota_{|z|>|w|}f(z,w)\in\mathbb{C}[[z,z^{-1},w,w^{-1}]]$ denotes the
power series expansion of a rational function $f$ in the region
$|z|>|w|$. We usually omit the symbol $\iota_{|z|>|w|}$ and just
write $(z-w)^{-1}$ to mean the expansion in the region $|z|>|w|$,
and write $-(w-z)^{-1}$ to mean the expansion in $|w|>|z|$. It is
easy to check that $a(w)\circ_n b(w)$ above is a well-defined
element of $QO(V)$.

The nonnegative circle products are connected through the {\it
operator product expansion} (OPE) formula.
For $a,b\in QO(V)$, we have \begin{equation} \label{opeformula} a(z)b(w)=\sum_{n\geq 0}a(w)\circ_n
b(w)~(z-w)^{-n-1}+:a(z)b(w):,\end{equation} which is often written as
$a(z)b(w)\sim\sum_{n\geq 0}a(w)\circ_n b(w)~(z-w)^{-n-1}$, where
$\sim$ means equal modulo the term $$
:a(z)b(w): \ =a(z)_-b(w)\ +\ (-1)^{|a||b|} b(w)a(z)_+.$$ Here
$a(z)_-=\sum_{n<0}a(n)z^{-n-1}$ and $a(z)_+=\sum_{n\geq
0}a(n)z^{-n-1}$. Note that $:a(w)b(w):$ is a well-defined element of
$QO(V)$. It is called the {\it Wick product} of $a$ and $b$, and it
coincides with $a\circ_{-1}b$. The other negative circle products
are related to this by
$$ n!~a(z)\circ_{-n-1}b(z)=\ :(\partial^n a(z))b(z):\ ,$$
where $\partial$ denotes the formal differentiation operator
$\frac{d}{dz}$. For $a_1(z),\dots ,a_k(z)\in QO(V)$, the $k$-fold
iterated Wick product is defined to be
$$ :a_1(z)a_2(z)\cdots a_k(z): = :a_1(z)b(z):,$$
where $b(z)= :a_2(z)\cdots a_k(z):$. We often omit the formal variable $z$ when no confusion will arise.

The set $QO(V)$ is a nonassociative algebra with the operations
$\circ_n$ and a unit $1$. We have $1\circ_n a=\delta_{n,-1}a$ for
all $n$, and $a\circ_n 1=\delta_{n,-1}a$ for $n\geq -1$. A linear subspace $\cA\subset QO(V)$ containing 1 which is closed under the circle products will be called a {\it quantum operator algebra} (QOA).
In particular $\cA$ is closed under $\partial$
since $\partial a=a\circ_{-2}1$. Many formal algebraic
notions are immediately clear: a homomorphism is just a linear
map that sends $1$ to $1$ and preserves all circle products; a module over $\cA$ is a
vector space $M$ equipped with a homomorphism $\cA\rightarrow
QO(M)$, etc. A subset $S=\{a_i|\ i\in I\}$ of $\cA$ is said to generate $\cA$ if any element $a\in\cA$ can be written as a linear
combination of nonassociative words in the letters $a_i$, $\circ_n$, for
$i\in I$ and $n\in\mathbb{Z}$. We say that $S$ {\it strongly generates} $\cA$ if any $a\in\cA$ can be written as a linear combination of words in the letters $a_i$, $\circ_n$ for $n<0$. Equivalently, $\cA$ is spanned by the collection $\{ :\partial^{k_1} a_{i_1}(z)\cdots \partial^{k_m} a_{i_m}(z):| ~i_1,\dots,i_m \in I,~ k_1,\dots,k_m \geq 0\}$.

We say that $a,b\in QO(V)$ {\it quantum commute} if $(z-w)^N
[a(z),b(w)]=0$ for some $N\geq 0$. Here $[,]$ denotes the super bracket. This condition implies that $a\circ_n b = 0$ for $n\geq N$, so (\ref{opeformula}) becomes a finite sum. If $N$ can be chosen to be $0$, we say that $a,b$ commute. A commutative quantum operator algebra (CQOA) is a QOA whose elements pairwise quantum commute. Finally, the notion of a CQOA is equivalent to the notion of a vertex algebra. Every CQOA $\cA$ is itself a faithful $\cA$-module, called the {\it left regular
module}. Define
$$\rho:\cA\rightarrow QO(\cA),\ \ \ \ a\mapsto\hat a,\ \ \ \ \hat
a(\zeta)b=\sum_{n\in\mathbb{Z}} (a\circ_n b)~\zeta^{-n-1}.$$ Then $\rho$ is an injective QOA homomorphism,
and the quadruple of structures $(\cA,\rho,1,\partial)$ is a vertex
algebra in the sense of \cite{FLM}. Conversely, if $(V,Y,{\bf 1},D)$ is
a vertex algebra, the collection $Y(V)\subset QO(V)$ is a
CQOA. {\it We will refer to a CQOA simply as a
vertex algebra throughout the rest of this paper}.

\section{Recollections on chiral equivariant cohomology}
We briefly recall the construction of chiral equivariant cohomology, following the notation in \cite{LLI}\cite{LLSI}\cite{LLSII}. A {\it differential vertex algebra} (DVA) is a degree graded vertex algebra $\cA^*=\oplus_{p\in\mathbb{Z}}\cA^p$ equipped with a vertex algebra derivation $d$ of degree 1 such that $d^2=0$. A DVA will be called {\it degree-weight graded} if it has an additional $\mathbb{Z}_{\geq 0}$-grading by weight, which is compatible with the degree in the sense that $\cA^p=\oplus_{n\geq0}\cA^p[n]$. There is an auxiliary structure on a DVA which is analogous to the structure of a $G^*$-algebra in \cite{GS}. Associated to $\gg$ is a Lie superalgebra $\gs\gg$ defined to be the semidirect product $\gg\triangleright \gg^{-1}$, where $\gg^{-1}$ is a copy of the adjoint module in degree $-1$. The bracket in $\gs\gg$ is given by $[(\xi,\eta),(x,y)]=([\xi,x],[\xi,y]-[x,\eta])$, and $\gs\gg$ is equipped with a differential $d:(\xi,\eta)\mapsto(\eta,0)$. This differential extends to the loop algebra $\gs\gg[t,t^{-1}]$, and gives rise to a vertex algebra derivation on the corresponding current algebra $O(\gs\gg):= O(\gs\gg,0)$. Here $0$ denotes the zero bilinear form on $\gs\gg$. An $O(\gs\gg)$-algebra is a degree-weight graded DVA $\cA$ equipped with a DVA homomorphism $\rho:O(\gs\gg)\rightarrow \cA$, which we denote by $(\xi,\eta)\mapsto L_{\xi} + \iota_{\eta}$.

In \cite{LLI}, the chiral equivariant cohomology functor was defined on the category of $O(\gs\gg)$-algebras, and in \cite{LLSI} this functor was extended to a larger class of spaces which carry only a representation of the Lie subalgebra $\gs\gg[t]$ of $\gs\gg[t,t^{-1}]$. An $\gs\gg[t]$-module is a degree-weight graded complex $(\cA,d_\cA)$ equipped with a Lie algebra homomorphism $\rho:\gs\gg[t]\ra End~\cA$, which we denote by $(\xi,\eta) t^n\mapsto L_{\xi}(n) + \iota_{\xi}(n)$, $n\geq 0$. We also require that for all $x\in\gs\gg[t]$ we have $\rho(dx)=[d_\cA,\rho(x)]$, and $\rho(x)$ has degree 0 whenever $x$ is even in $\gs\gg[t]$, and degree -1 whenever $x$ is odd, and has weight $-n$ if $x\in\gs\gg t^n$. Finally, we require $\cA$ to admit a compatible action $\hat{\rho}:G\ra Aut(\cA)$ of $G$ satisfying:
\begin{equation} \label{gsI} \frac{d}{dt}\hat{\rho}(exp(t\xi))|_{t=0} = L_{\xi}(0),\end{equation}
\begin{equation} \label{gsII} \hat{\rho} (g) L_{\xi}(n) \hat{\rho}(g^{-1}) = L_{Ad(g)(\xi)}(n),\end{equation}
\begin{equation} \label{gsIII} \hat{\rho}(g) \iota_{\xi}(n)\hat{\rho}(g^{-1}) =\iota_{Ad(g)(\xi)}(n),\end{equation}
\begin{equation} \label{gsIV} \hat{\rho}(g) d \hat{\rho} (g^{-1}) = d,\end{equation}
for all $\xi\in\gg$, $g\in G$, and $n\geq 0$. These conditions are analogous to Equations (2.23)-(2.26) of \cite{GS}. In order for (\ref{gsI}) to make sense, we must be able to differentiate along appropriate curves in $\cA$, which is the case in our main examples. Given an $\gs\gg[t]$-module $(\cA,d)$, we define the chiral horizontal, invariant and basic subspaces of $\cA$ to be respectively
$$
\cA_{hor}=\{a\in\cA|\rho(x)a=0,\ \forall x\in\gg^{-1}[t]\}, $$ $$
\cA_{inv}=\{a\in\cA|\rho(x)a=0,~\forall x\in\gg[t],~\hat{\rho}(g)(a) = a,~\forall g\in G\},$$ $$
\cA_{bas}=\cA_{hor}\cap\cA_{inv}.$$ Both $\cA_{inv}$ and $\cA_{bas}$ are subcomplexes of $\cA$, but $\cA_{hor}$ is not a subcomplex of $\cA$ in general.

An $O(\gs\gg)$-algebra which plays an important role in our theory is the semi-infinite Weil complex $\cW(\gg)$, which is just the $bc\beta\gamma$ system $\cE(\gg)\otimes\cS(\gg)$. The $bc$ and $\beta\gamma$ systems were introduced by Friedan-Martinec-Skenker in \cite{FMS}; in this notation, the $bc$ system $\cE(\gg)$ is the vertex algebra with odd generators $b^{\xi}, c^{\xi'}$, which are linear in $\xi \in\gg$ and $\xi'\in \gg^*$, and satisfy the OPE relations $$b^{\xi}(z) c^{\xi'}(w)\sim\langle \xi',\xi\rangle (z-w)^{-1},\ \ \ \ \ \ c^{\xi'}(z)b^{\xi}(w)\sim \langle \xi',\xi\rangle (z-w)^{-1},$$ $$b^{\xi}(z)b^{\eta}(w)\sim 0,\ \ \ \ \ \ c^{\xi'}(z) c^{\eta'}(w)\sim 0.$$ Here $\langle,\rangle$ denotes the natural pairing between $\gg^*$ and $\gg$. Similarly, the $\beta\gamma$ system $\cS(\gg)$ is the vertex algebra with even generators $\beta^{\xi}, \gamma^{\xi'}$, which are linear in $\xi \in\gg$ and $\xi'\in \gg^*$, and satisfy $$\beta^{\xi}(z) \gamma^{\xi'}(w)\sim\langle \xi',\xi\rangle (z-w)^{-1},\ \ \ \ \ \ \gamma^{\xi'}(z)\beta^{\xi}(w)\sim -\langle \xi',\xi\rangle (z-w)^{-1},$$ $$\beta^{\xi}(z)\beta^{\eta}(w)\sim 0,\ \ \ \ \ \ \gamma^{\xi'}(z)\gamma^{\eta'}(w)\sim 0.$$ 
$\cW(\gg)$ possesses a Virasoro element $\omega_{\cW}$ of central charge zero, given by 
\begin{equation}\label{weilvir}
\omega_\cW=\omega_\cE+\omega_\cS,\ \ \ \ \ \ \ \ \ 
\omega_\cE=-\sum_{i=1}^n :b^{\xi_i}\partial c^{\xi_i'}:, \ \ \ \ \ \ \ \ \ 
\omega_\cS= \sum_{i=1}^n :\beta^{\xi_i}\partial\gamma^{\xi_i'}:.
\end{equation} Here $\{\xi_1,\dots,\xi_n\}$ is a basis for $\gg$ and $\{\xi'_1,\dots,\xi'_n\}$ denotes the corresponding dual basis  for $\gg^*$. The generators $\beta^{\xi},\gamma^{\xi'}, b^{\xi}, c^{\xi'}$ are primary of weights $1,0,1,0$ with respect to $\omega_{\cW}$. There is an additional grading by degree, in which $\beta^{\xi},\gamma^{\xi'}, b^{\xi}, c^{\xi'}$ (and their respective derivatives) have degrees $-2,2,-1,1$. Note that the weight zero component is isomorphic to the classical Weil algebra $W(\gg)=\bigwedge(\gg^*)\otimes Sym(\gg^*)$, where the degree 1 and degree 2 generators $c^{\xi'}, \gamma^{\xi'}$ play the role of connection 1-forms and curvature 2-forms, respectively.

Next, we recall the $O(\gs\gg)$-algebra structure on $\cW(\gg)$. Define vertex operators \begin{equation} \label{defoftheta}
\Theta_\cW^\xi=\Theta_\cE^\xi+\Theta_\cS^\xi,\ \ \ \ 
\Theta_\cE^\xi= \sum_{i=1}^n :b^{[\xi,\xi_i]}c^{\xi_i'}:,\ \ \ \ 
\Theta_\cS^\xi=- \sum_{i=1}^n :\beta^{[\xi,\xi_i]}\gamma^{\xi_i'}:,\end{equation}
\begin{equation} \label{defofdiff}
 D=J+K,\ \ \ \ J= \sum_{i=1}^n :(\Theta^{\xi_i}_{\cS} + \frac{1}{2}\Theta^{\xi_i}_{\cE})c^{\xi'_i}:,\ \ \ \ K= \sum_{i=1}^n:\gamma^{\xi_i'}b^{\xi_i}:.
\end{equation}

The Fourier modes $J(0)$, $K(0)$, and $D(0)$ are called the BRST, chiral Koszul, and chiral Weil differentials, respectively. They satisfy $J(0)^2=K(0)^2=D(0)^2=[K(0),J(0)]=0$. Moreover, we have
$$\Theta^{\xi}_{\cW}(z)  \Theta^{\eta}_{\cW}(w) \sim \Theta^{[\xi,\eta]}_{\cW}(w)(z-w)^{-1},$$ $$ \Theta^{\xi}_{\cW} (z) b^{\eta}(w) \sim b^{[\xi,\eta]}(w)(z-w)^{-1},$$ $$[D(0),b^{\xi}(z)] = \Theta^{\xi}_{\cW}(z),\ \ \ \ \  [D(0), \Theta^{\xi}_{\cW}(z)] = 0.$$ In other words, the map $O(\gs\gg) \rightarrow \cW(\gg)$ sending $(\xi,\eta)(z) \mapsto \Theta^{\xi}_{\cW}(z) + b^{\eta}(z)$ and sending $d\mapsto [D(0),-]$, is a homomorphism of DVAs.

The horizontal subalgebra $\cW(\gg)_{hor}$ is the vertex subalgebra generated by $\beta^{\xi},\gamma^{\xi'},b^{\xi}$. The basic subalgebra $\cW(\gg)_{bas}=\cW(\gg)_{hor}^{\gg[t]}$ is easily seen to be a subcomplex under both differentials $K(0)$ and $J(0)$.

\begin{defn}For any $\gs\gg[t]$-module $(\cA,d_{\cA})$, we define its chiral basic cohomology $\H^*_{bas}(\cA)$ to be $ H^*(\cA_{bas},d_{\cA})$. We define its chiral equivariant cohomology $\H^*_{G}(\cA)$ to be $\H^*_{bas}(\cW(\gg)\otimes\cA)$. \end{defn}

In this paper, we are interested in the case when $\cA$ is the trivial $\gs\gg[t]$-module $\mathbb{C}$. We refer to $\H^*_G(\mathbb{C})$ as the chiral point algebra; it plays the role of $H^*_G(pt) = Sym(\gg^*)^G$ in the classical theory, and $\H^*_G(\mathbb{C})[0]\cong H^*_G(pt)$. 

Any connected, compact group $G$ can be written as a quotient $$(G_1\times \cdots \times G_r \times T)/\Gamma,$$ where the $G_i$ are simple, $T$ is a torus, and $\Gamma$ is finite. Then $\H^*_G(\mathbb{C}) = \H^*_{G_1}(\mathbb{C}) \otimes \cdots \otimes \H^*_{G_r}(\mathbb{C}) \otimes \H^*_T(\mathbb{C})$. By Theorem 6.1 of \cite{LLI}, $\H^*_T(\mathbb{C})$ is the free abelian vertex algebra with generators $\gamma^{\xi'_1},\dots,\gamma^{\xi'_n}$, where $\{\xi_1,\dots,\xi_n\}$ is a basis for the Lie algebra of $T$. In other words, $\H^*_T(\mathbb{C})$ is the polynomial algebra generated by $\partial^k\gamma^{\xi'_i}$ for $k\geq 0$ and $i=1,\dots,n$. So we may assume without loss of generality that $G$ is simple. 

There is a \lq\lq classical sector" of $\H^*_G(\mathbb{C})$, which is the abelian subalgebra generated by the weight zero component $\H^*_G(\mathbb{C})[0] = Sym(\gg^*)^G$. Unlike the case where $G$ is abelian, $\H^*_G(\mathbb{C})$ contains additional elements that have no classical analogues. The most notable feature of $\H^*_G(\mathbb{C})$ is a conformal structure of central charge zero. Let $\{\xi_1,\dots,\xi_n\}$ be an orthonormal basis for $\gg$ relative to the Killing form. The Virasoro element $\L$ is represented by \begin{equation}\label{defofvirasoro} L = \omega_{\cS} -L_{\cS} + C\end{equation} where $L_{\cS} = - \sum_{i=1}^n :\Theta^{\xi_i}_{\cS}  \Theta^{\xi_i}_{\cS}:$ and $C = \sum_{i,j=1}^n :b^{\xi_i} b^{\xi_j} \gamma^{ad^*(\xi_i)(\xi'_j)}: = :(K(0)\Theta^{\xi_i}_{\cS}) b^{\xi_i}:$. The term $ \omega_{\cS} -L_{\cS}$ lies in $\cW(\gg)_{bas}$ and satisfies the Virasoro OPE with central charge zero, but it is not $D(0)$-closed. The purpose of the term $C$ is to correct this flaw, and $L$ still satisfies $L(z) L(w) \sim 2L(w) (z-w)^{-2} + \partial L(w) (z-w)^{-1}$. By Corollary 7.17 of \cite{LLI}, $L$ represents a nontrivial class $\L$ in $\H^*_G(\mathbb{C})$, and by Corollary 4.18 of \cite{LLSI}, $\L$ is a conformal structure on $\H^*_G(\mathbb{C})$. The proof is very simple. First, the Virasoro element $\omega_{\cW} = \omega_{\cS} + \omega_{\cE}$ is a conformal structure on $\cW(\gg)$; in particular, $\omega_{\cW}\circ_0$ acts by $\partial$ and $\omega_{\cW} \circ_1 a = m a$ for all $a\in \cW(\gg)[m]$. Even though $\omega_{\cW}$ does not lie in $\cW(\gg)_{bas}$, the operators $\omega_{\cW} \circ_k$ act on $\cW(\gg)_{bas}$ for all $k\geq 0$. For all $a\in \cW(\gg)_{bas}$ and $k\geq 0$, we have \begin{equation} \label{defofh}(L - \omega_{\cW})\circ_k a = D(0) (H\circ_k a),\ \ \ \ \ \ \ H = \sum_{i=1}^n :\Theta^{\xi_i}_{\cS} b^{\xi_i}:\end{equation} Even though $H$ does not lie in $\cW(\gg)_{bas}$, $H \circ_k$ preserves $\cW(\gg)_{bas}$ for all $k\geq 0$. Hence the operators $L\circ_k$ and $\omega_{\cW}\circ_k$ on $\cW(\gg)_{bas}$ agree up to coboundary. 

The Virasoro element $\L$ is nonclassical and has no analogue in $H^*_G(pt)$, since $H^0_G(pt) = \mathbb{C}$. The main purpose of this paper is to construct in a uniform way an infinite family of new, nonclassical elements of $\H^*_G(\mathbb{C})$ in the case $G=SU(2)$. In particular, we construct an injective linear map $\Psi:U(\gs\gl_2)\rightarrow \H^*_{SU(2)}(\mathbb{C})$, for which $\Psi(1) = \L$. The image of $\Psi$, together with the generator of $H^*_{SU(2)}(pt)$, is conjectured to be a strong generating set for $\H^*_{SU(2)}(\mathbb{C})$. We also compute some relations among these generators, which gives some glimpse of the rich algebraic structure possessed by the chiral point algebra.

\section{A peculiar topological structure on $\cW(\gg)$}

The notion of a topological vertex algebra (TVA) was introduced by Lian-Zuckerman in \cite{LZII}. It is an abstraction based on examples from physics. A TVA is a vertex algebra $\cA$ equipped with four distinguished vertex operators $L,F,J,G$, where $L$ is a Virasoro element with central charge zero, $F$ is an even current which is conformal weight one quasi-primary (with respect to $L$), $J$ is an odd conformal weight one primary, and $G$ is an odd conformal weight two primary, such that
\begin{equation}\label{deftop}J(z)J(w)\sim 0,\ \ \ \ \ \ G(z)G(w)\sim 0,\ \ \ \ \ \ J(0)G=L,\ \ \ \ \ \  F(0)J=J,\ \ \ \ \ \ F(0)G=-G.\end{equation}
We do not require that $L$ is a conformal structure on $\cA$, and in particular, $L_0 = L\circ_1$ need not act diagonalizably on $\cA$. 

There is a well-known TVA structure on the $bc\beta\gamma$ system $\cE(V)\otimes \cS(V)$ attached to an $n$-dimensional vector space $V$ with basis $\{x_1,\dots,x_n\}$. Define elements 
$$L = \sum_{i=1}^n \big(:\beta^{x_i} \partial \gamma^{x'_i} :-  :b^{x_i} \partial c^{x'_i}: \big), \ \ \ \ F = - \sum_{i=1}^n :b^{x_i} c^{x'_i}:, \ \ \ J= \sum_{i=1}^n :c^{x'_i} \beta^{x_i}:, \ \ \ \ G= -\sum_{i=1}^n :b^{x_i} \partial\gamma^{x'_i}:. $$ It is easy to check that (\ref{deftop}) is satisfied. In particular, for any Lie algebra $\gg$, $\cW(\gg)$ possesses this TVA structure.

The same vertex algebra $\cW(\gg)$ supports another, more subtle TVA structure, which depends on the $O(\gs\gg)$ structure. As usual, fix an orthonormal basis $\{\xi_1,\dots, \xi_n\}$ for $\gg$ relative to the Killing form.
\begin{thm} Define \begin{equation}\label{newtopi} L =  -\sum_{i=1}^n :\Theta^{\xi_i}_{\cS} \Theta^{\xi_i}_{\cS}: - \sum_{i=1}^n :\Theta^{\xi_i}_{\cW} \Theta^{\xi_i}_{\cW}: - \sum_{i=1}^n :b^{\xi_i} \partial c^{\xi'_i}:,\end{equation} \begin{equation} \label{newtopii} J =\sum_{i=1}^n : (\Theta^{\xi_i}_{\cS} + \frac{1}{2}\Theta^{\xi_i}_{\cE}) c^{\xi'_i}:,\ \ \ \ \ \ G  = -\sum_{i=1}^n : (2\Theta^{\xi_i}_{\cS} + \Theta^{\xi_i}_{\cE}) b^{\xi_i}:,\ \ \ \ \ \ F = - \sum_{i=1}^n :b^{\xi_i} c^{\xi'_i}:.\end{equation} These elements commute with $\Theta^{\xi}_{\cW}$ for all $\xi\in\gg$ and satisfy (\ref{deftop}), and hence define a TVA structure inside $\cW(\gg)^{\gg[t]}$. \end{thm}

\begin{proof} This is a straightforward OPE calculation. \end{proof}

Note that these vertex operators do not lie in $\cW(\gg)_{bas}$, since they depend on $c^{\xi'_i}$. However, the nonnegative modes $F(k)$, $J(k)$, $G(k)$, $L(k)$ for $k\geq 0$ act on $\cW(\gg)_{bas}$. This follows readily from the OPEs
$$J(z) b^{\xi}(w)\sim \Theta^{\xi}_{\cW}(w)(z-w)^{-1},$$
$$G(z) b^{\xi}(w) \sim -\sum_{i=1}^n b^{[\xi_i,\xi]} b^{\xi_i}(w)(z-w)^{-1},$$
$$F(z) b^{\xi}(w)\sim -b^{\xi}(w) (z-w)^{-1}.$$ This structure on $\cW(\gg)_{bas}$ will be useful later in the study of $\H^*_G(\mathbb{C})$ for $G=SU(2)$.

\section{Graded and filtered structures}
Let $\cR$ be the category of vertex algebras $\cA$ equipped with a $\mathbb{Z}_{\geq 0}$-filtration
\begin{equation} \cA_{(0)}\subset\cA_{(1)}\subset\cA_{(2)}\subset \cdots,\ \ \ \cA = \bigcup_{k\geq 0}
\cA_{(k)}\end{equation} such that $\cA_{(0)} = \mathbb{C}$, and for all
$a\in \cA_{(k)}$, $b\in\cA_{(l)}$, we have
\begin{equation} \label{goodi} a\circ_n b\in\cA_{(k+l)},\ \ \ \text{for}\
n<0,\end{equation}
\begin{equation} \label{goodii} a\circ_n b\in\cA_{(k+l-1)},\ \ \ \text{for}\
n\geq 0.\end{equation}
Elements $a(z)\in\cA_{(d)}\setminus \cA_{(d-1)}$ are said to have degree $d$.

Filtrations on vertex algebras satisfying (\ref{goodi})-(\ref{goodii})~were introduced in \cite{LiII}, and are known as {\it good increasing filtrations}. If $\cA$ possesses such a filtration, the associated graded object $gr(\cA) = \bigoplus_{k>0}\cA_{(k)}/\cA_{(k-1)}$ is a
$\mathbb{Z}_{\geq 0}$-graded associative, supercommutative algebra with a
unit $1$ under a product induced by the Wick product on $\cA$. In general, there is no natural linear map $\cA\rightarrow gr (\cA)$, but for each $r\geq 1$ we have the projection \begin{equation} \phi_r: \cA_{(r)} \ra \cA_{(r)}/\cA_{(r-1)}\subset gr(\cA).\end{equation} 
Moreover, $gr(\cA)$ has a derivation $\partial$ of degree zero
(induced by the operator $\partial = \frac{d}{dz}$ on $\cA$), and
for each $a\in\cA_{(d)}$ and $n\geq 0$, the operator $a\circ_n$ on $\cA$
induces a derivation of degree $d-k$ on $gr(\cA)$, which we denote by $a(n)$. Here $$k  = sup \{ j\geq 1|~ \cA_{(r)}\circ_n \cA_{(s)}\subset \cA_{(r+s-j)}~\forall r,s,n\geq 0\},$$ as in \cite{LL}. Finally, these derivations give $gr(\cA)$ the structure of a vertex Poisson algebra.

The assignment $\cA\mapsto gr(\cA)$ is a functor from $\cR$ to the category of $\mathbb{Z}_{\geq 0}$-graded supercommutative rings with a differential $\partial$ of degree 0, which we will call $\partial$-rings. A $\partial$-ring is the same thing as an {\it abelian} vertex algebra, that is, a vertex algebra $\cV$ in which $[a(z),b(w)] = 0$ for all $a,b\in\cV$. A $\partial$-ring $A$ is said to be generated by a subset $\{a_i|~i\in I\}$ if $\{\partial^k a_i|~i\in I, k\geq 0\}$ generates $A$ as a graded ring. The key feature of $\cR$ is the following reconstruction property \cite{LL}:

\begin{lemma}\label{reconlem}Let $\cA$ be a vertex algebra in $\cR$ and let $\{a_i|~i\in I\}$ be a set of generators for $gr(\cA)$ as a $\partial$-ring, where $a_i$ is homogeneous of degree $d_i$. If $a_i(z)\in\cA_{(d_i)}$ are vertex operators such that $\phi_{d_i}(a_i(z)) = a_i$, then $\cA$ is strongly generated as a vertex algebra by $\{a_i(z)|~i\in I\}$.\end{lemma}

For any Lie algebra $\gg$, the semi-infinite Weil algebra $\cW(\gg)$ admits a good increasing filtration \begin{equation} \label{filtw} \cW(\gg)_{(0)}\subset \cW(\gg)_{(1)}\subset \cdots,\ \ \ \ \ \ \ \ \cW(\gg) = \bigcup_{k\geq 0} \cW(\gg)_{(k)},\end{equation} where $\cW(\gg)_{(k)}$ is defined to be the vector space spanned by iterated Wick products of the generators $b^{\xi}, c^{\xi'}, \beta^{\xi}, \gamma^{\xi'}$ and their derivatives, of length at most $k$. We say that elements of $\cW(\gg)_{(k)}\setminus \cW(\gg)_{(k-1)}$ have {\it polynomial degree $k$}. This filtration is $\gs\gg[t]$-invariant, and we have an isomorphism of supercommutative rings \begin{equation}\label{assgrad} gr(\cW(\gg)) \cong Sym(\bigoplus_{k\geq 0} (V_k \oplus V^*_k)) \bigotimes \bigwedge (\bigoplus_{k\geq 0} (U_k \oplus U^*_k)).\end{equation} Here $V_k$, $U_k$ are copies of $\gg$, and $V^*_k$, $U^*_k$ are copies of $\gg^*$. The generators of $gr(\cW(\gg))$ are $\beta^{\xi_i}_k$, $\gamma^{\xi'_i}_k$, $b^{\xi_i}_k$, and $c^{\xi'_i}_k$, which correspond to the vertex operators $\partial^k \beta^{\xi_i}$, $\partial^k \gamma^{\xi'_i}$, $\partial^k b^{\xi_i}$, and $\partial^k c^{\xi'_i}$, respectively for $k\geq 0$. Since $\cW(\gg)$ has a basis consisting of iterated Wick products of the generators and their derivatives, $\cW(\gg)\cong gr(\cW(\gg))$ as vector spaces, although not canonically. Finally, the filtration (\ref{filtw}) is inherited by any subalgebra of $\cW(\gg)$, such as $\cW(\gg)_{hor}$ and $\cW(\gg)_{bas}$. Note that \begin{equation}\label{parity} K(0)(\cW(\gg)_{(k)}) \subset \cW(\gg)_{(k)},\ \ \ \ \ \ \ \ J(0)(\cW(\gg)_{(k)}) \subset \cW(\gg)_{(k+1)}. \end{equation}  Clearly $K(0)$ preserves the parity of polynomial degree and $J(0)$ reverses it.

Recall that the horizontal subalgebra $\cW(\gg)_{hor}\subset \cW(\gg)$  is generated by $\beta^{\xi},\gamma^{\xi'}, b^{\xi}$. Clearly $\cW(\gg)_{hor}$ (but not $\cW(\gg)$) has a $\mathbb{Z}_{\geq 0}$ grading by $b$-number, which is the eigenvalue of the Fourier mode $-F(0)$, where $F$ is given by (\ref{newtopii}). This grading is inherited by $\cW(\gg)_{bas}$, since the action of $\Theta^{\xi}_{\cW}$ on $\cW(\gg)_{hor}$ preserves $b$-number, for all $\xi\in\gg$. Let us introduce the notations $\cW(\gg)_{hor}^{(k)}$ and $\cW(\gg)_{bas}^{(k)}$ to denote the subspaces of $b$-number $k$; in this notation, $\cS(\gg)^{\gg[t]}=\cW(\gg)_{bas}^{(0)}$. The parity of the $b$-number and cohomological degree gradations coincide, but since $K(0)$ and $J(0)$ raise and lower $b$-number by $1$, respectively, $\H^*_G(\mathbb{C})$ is not graded by $b$-number. 

However, the $b$-number gradation does induce a filtration on $\H^*_G(\mathbb{C})$. Any class in $\H^{2k}_G(\mathbb{C})$ has a representative of the form $\omega = \sum_{i\geq 0} \omega^{(2i)}$, where $\omega^{(2i)}\in \cW(\gg)_{bas}^{(2i)}$, and only finitely many terms are nonzero. Consequently, for each $k\in\mathbb{Z}$, $\H^{2k}_G(\mathbb{C})$ admits a decreasing filtration \begin{equation}\label{bnumbfiltrationi} \H^{2k}_G(\mathbb{C})_{(0)} \supset  \H^{2k}_G(\mathbb{C})_{(2)} \supset   \H^{2k}_G(\mathbb{C})_{(4)} \supset \cdots, \end{equation} where $\H^{2k}_G(\mathbb{C})_{(2r)}$ consists of classes admitting a representative of the form $\omega = \sum_{i\geq r} \omega^{(2i)}$. There is a similar filtration on $\H^{2k+1}_G(\mathbb{C})$ of the form
$$ \H^{2k+1}_G(\mathbb{C})_{(1)} \supset   \H^{2k+1}_G(\mathbb{C})_{(3)}  \supset \H^{2k+1}_G(\mathbb{C})_{(5)} \supset \cdots. $$
Note that if $\omega = \sum_{i\geq l} \omega^{(2i)}$ represents a class in $\H^{2k}_G(\mathbb{C})$, we must have $J(0)(\omega^{(2l)}) = 0$, since $K(0)$ and $J(0)$ raise and lower $b$-number by 1, respectively. Hence the lowest term $\omega^{(2l)}$ represents a class in $H^{2k}(\cW(\gg)_{bas},J(0))$, and we obtain an injective linear map \begin{equation}\label{injbj} \H^{2k}_G(\mathbb{C})_{(2l)}/ \H^{2k}_G(\mathbb{C})_{(2l+2)} \rightarrow H^{2k}(\cW(\gg)_{bas},J(0)),\end{equation} whose image is homogeneous $b$-number $2l$. In odd degree $2k+1$, there is a similar injective map $$\H^{2k+1}_G(\mathbb{C})_{(2l+1)} / \H^{2k+1}_G(\mathbb{C})_{(2l+3)}\rightarrow H^{2k+1}(\cW(\gg)_{bas},J(0)),$$ whose image is homogeneous of $b$-number $2l+1$. 

\section{The basic complex and the vertex algebra commutant problem}

The main technical difficulty that arises in studying $\H^*_G(\mathbb{C})$ is that the basic complex is a commutant subalgebra of $\cW(\gg)$, and vertex algebra commutants are generally not well understood. Given a vertex algebra $\cV$ and a subalgebra $\cA$, recall that the commutant $Com(\cA,\cV)$ is defined to be $$\{v\in \cV| [a(z),v(w)]=0,\forall a\in\cA\}.$$ If $\cA$ is a homomorphic image of a current algebra $O(\gg,B)$ of some Lie algebra $\gg$, we have $Com(\cA,\cV) = \cV^{\gg[t]}$. In this notation, $\cW(\gg)_{bas} = Com(\cO,\cW(\gg))$, where $\cO$ is the copy of $O(\gs\gg)$ generated by $\Theta^{\xi}_{\cW}, b^{\xi}$ for $\xi \in\gg$. It is difficult to study this algebra because $O(\gs\gg)$ does not act completely reducibly on $\cW(\gg)$. Moreover, it seems impossible to describe $\H^*_G(\mathbb{C})$ without first giving a reasonable description of $\cW(\gg)_{bas}$. Even in the simplest case $G=SU(2)$, we are unable to describe $\cW(\gg)_{bas}$ completely. However, in this case we will find a generating set for the subalgebra $\cW(\gg)_{bas}^{(0)} = \cS(\gg)^{\gg[t]}$, as well as for the subspace $\cW(\gg)_{bas}^{(1)}$, which is a module over $\cS(\gg)^{\gg[t]}$. This will be sufficient to construct the injective linear map $\Psi: U(\gs\gl_2)\rightarrow \H^*_{SU(2)}(\mathbb{C})$.

As we shall see, $\cS(\gg)^{\gg[t]}$ plays an important role in the structure of the chiral point algebra, but unfortunately we lack the tools to describe $\cS(\gg)^{\gg[t]}$ for any simple $\gg$ other than $\gs\gl_2$. However, since $\cS(\gg)^{\gg[t]}$ coincides with $Com(\Theta_{\cS},\cS(\gg))$, where $\Theta_{\cS}$ is generated by $\{\Theta^{\xi}_{\cS}|~\xi\in\gg\}$, we have the following simple characterization of this algebra: 

\begin{lemma}\label{bnumbzero} $\cS(\gg)^{\gg[t]} = \cS(\gg)\cap Ker(J(0))$. In particular, we have a vertex algebra homomorphism $\cS(\gg)^{\gg[t]}\rightarrow H^*(\cW(\gg)_{bas},J(0))$. \end{lemma}

\begin{proof} An OPE calculation shows that $J\circ_0\omega = \frac{1}{m!} \sum_{m\geq 0} :(\partial^m c^{\xi'_i})(\Theta^{\xi_i}_{\cS}\circ_m \omega):$ for any $\omega\in \cS(\gg)^{\gg[t]}$, so the claim is immediate. \end{proof}

$\cS(\gg)^{\gg[t]}$ possesses some additional features that were described in \cite{LL}. In terms of an orthonormal basis $\{\xi_1,\dots,\xi_n\}$ for $\gg$ relative to the Killing form, there is a Virasoro element of central charge zero, given by $$\omega_{\cS} - L_{\cS} = \sum_{i=1}^n :\beta^{\xi_i} \partial \gamma^{\xi'_i}:+ \sum_{i=1}^n :\Theta^{\xi_i}_{\cS} \Theta^{\xi_i}_{\cS}: .$$ There is also an action of the current algebra $O(\gs\gl_2,-\frac{n}{8}\kappa)$ with generators $X^x, X^h, X^y$, given by \begin{equation}\label{basicslii} X^h\mapsto \sum_{i=1}^n :\beta^{\xi_i} \gamma^{\xi'_i}:,\ \ \ \ X^x\mapsto \frac{1}{2} \sum_{i=1}^n:\gamma^{\xi'_i} \gamma^{\xi'_i}:,\ \ \ \ \ X^y\mapsto -\frac{1}{2} \sum_{i=1}^n:\beta^{\xi_i} \beta^{\xi_i}:.\end{equation} In fact, if $V$ is any $\gg$-module of dimension $m$ which admits a symmetric, invariant bilinear form, $\cS(V)^{\gg[t]}$ carries an action of $O(\gs\gl_2,-\frac{m}{8}\kappa)$ given by the same formula. Moreover, $O(\gs\gl_2,-\frac{m}{8}\kappa)$ has level $-\frac{m}{2}$ in the standard normalization. By Theorem 0.2.1 of \cite{GK}, it is a {\it simple} vertex algebra, so the map $O(\gs\gl_2,-\frac{m}{8}\kappa)\rightarrow \cS(V)^{\gg[t]}$ is injective.

In \cite{LL}, an approach to studying vertex algebras of the form $\cS(V)^{\gg[t]}$ using commutative algebra was introduced. There is a good increasing filtration \begin{equation}\label{filts}\cS(V)_{(0)}\subset \cS(V)_{(1)}\subset \cS(V)_{(2)} \subset \cdots,\ \ \ \ \ \cS(V) = \bigcup_{k\geq 0} \cS(V)_{(k)}, \end{equation} where $\cS(V)_{(k)}$ is spanned by iterated Wick products of $\beta^{x}, \gamma^{x'}$ and their derivatives, of length at most $k$. This is analogous to the Bernstein filtration on the Weyl algebra $\cD(V)$. This filtration is $\gg[t]$-invariant, so $\gg[t]$ acts on $gr(\cS(V))$ by derivations of degree zero, and there is an injective map of invariant spaces \begin{equation} \label{defofgamma}  gr(\cS(V)^{\gg[t]})\hookrightarrow gr(\cS(V))^{\gg[t]}.\end{equation} The latter space may be easier to describe, and if this map happens to be surjective we can reconstruct $\cS(V)^{\gg[t]}$ using Lemma \ref{reconlem}, since a set of generators for $gr(\cS(V)^{\gg[t]})$ as a differential algebra corresponds to a strong generating set for $\cS(V)^{\gg[t]}$ as a vertex algebra. 

The problem of describing $gr(\cS(V))^{\gg[t]}$ can be reinterpreted in the language of jet schemes. Let $X$ be an irreducible scheme of finite type over $\mathbb{C}$. For each integer $m\geq 0$, the jet scheme $J_m(X)$ is determined by its functor of points: for every $\mathbb{C}$-algebra $A$, we have a bijection
$$Hom (Spec (A), J_m(X)) \cong Hom (Spec (A[t]/\langle t^{m+1}\rangle ), X).$$ Thus the $\mathbb{C}$-valued points of $J_m(X)$ correspond to the $\mathbb{C}[t]/\langle t^{m+1}\rangle$-valued points of $X$. We define $J_{\infty}(X) = \lim_{\infty \leftarrow m} J_m(X)$, which is known as the infinite jet scheme, or space of arcs of $X$, and we denote by $\cO(J_{\infty}(X))$ the ring $\lim_{m\rightarrow \infty} \cO(J_m(X))$. It is a differential algebra with derivation $D$. More details about jet schemes and a list of references can be found in the Appendix.

For a finite-dimensional vector space $V$, there is an isomorphism of differential algebras $$gr(\cS(V))\cong \cO(J_{\infty}(V\oplus V^*))$$ which intertwines the differentials $D$ and $\partial$. If $V$ is a representation of a (connected) Lie group $G$ with Lie algebra $\gg$, there is an action of $\gg[t]$ on $\cO(J_{\infty}(V\oplus V^*))$, and this isomorphism intertwines the actions of $\gg[t]$ as well. Hence we have an isomorphism of invariant spaces $$gr(\cS(V))^{\gg[t]}\cong \cO(J_{\infty}(V\oplus V^*))^{\gg[t]}.$$

We are thus led to the problem of describing invariant rings of the form $\cO(J_{\infty}(V))^{\gg[t]}$, where $V$ is an arbitrary finite-dimensional $G$-representation. Note that $\cO(V)^G$ is a canonical subalgebra of $\cO(J_{\infty}(V))^{\gg[t]}$. Let $\langle \cO(V)^G\rangle$ denote the algebra generated by $\{D^i (f)|~f\in \cO(V)^G, ~i\geq 0\}$, which lies in $\cO(J_{\infty}(V))^{\gg[t]}$, since the latter is closed under $D$. The following result is important in our study of $\H^*_{SU(2)}(\mathbb{C})$, since it will allow us to give a complete description of $\cS(\gg)^{\gg[t]}$ and $\cW(\gg)_{bas}^{(1)}$, for $G = SU(2)$.

\begin{thm} \label{jet} Let $G = SU(2)$ and let $V$ be the adjoint representation of $G$. Then we have \begin{equation} \label{jettwo}\cO(J_{\infty}(V\oplus V))^{\gg[t]} = \langle \cO(V\oplus V)^G\rangle,\end{equation} \begin{equation} \label{jetthree} \cO(J_{\infty}(V\oplus V\oplus V))^{\gg[t]} = \langle \cO(V\oplus V \oplus V)^G\rangle.\end{equation} \end{thm}

We will develop some general techniques for studying invariant rings of this kind and prove this result in the Appendix.

\section{The case $G=SU(2)$}
In this section, we consider the simplest nontrivial case $G=SU(2)$. The complexified Lie algebra of $SU(2)$ is $\gs\gl_2$, and we work in the standard root basis $x,y,h$, with commutators $$[x,y]=h,\ \ \ \ \ \  [h,x]=2x,\ \ \ \ \ \ [h,y]=-2y.$$ For simplicity of notation, we denote $\cW(\gs\gl_2)$ and $\cS(\gs\gl_2)$ by $\cW$ and $\cS$, respectively. We need the following theorem of Weyl which describes polynomial invariants for the adjoint representation of $\gs\gl_2$.

\begin{thm}\label{weyl}
For $n\geq 0$, let $V_n$ be a copy of the adjoint representation of $\gs\gl_2$, with basis $\{a^h_n,a^x_n,a^y_n\}$. Then $Sym\big(\bigoplus_{n=0}^{\infty}V_n\big)^{\gs\gl_2}$ is generated by
\begin{equation}\label{quadgen} q_{ij}= a^h_i a^h_j + 2a^x_i a^y_j + 2a^x_j a^y_i, \ \ \ \ \ \ \ c_{klm}= \left| \begin{array}{lll}
a^h_k & a^x_k & a^y_k \\ a^h_l & a^x_l & a^y_l \\ a^h_m & a^x_m &
a^y_m \end{array}\right|,\ \ \ \ \ \ \ i,j\geq 0,\ \ \ \ \ \ \ 0\leq k<l<m. \end{equation}

The ideal of relations among the variables $q_{ij}$ and $c_{klm}$
is generated by polynomials of the following two types: \begin{equation}\label{firstrel}
q_{ij}c_{klm}-q_{kj}c_{ilm}+q_{lj}c_{kim}-q_{mj}c_{kli}, \end{equation}
\begin{equation}\label{secrel} c_{ijk}c_{lmn}+\frac{1}{4}\left| \begin{array}{lll} q_{il} & q_{im} &
q_{in} \\q_{jl} & q_{jm} & q_{jn} \\ q_{kl} & q_{km} & q_{kn}\end{array}\right| .
\end{equation} Finally, the ideal of relations among the quadratics  $q_{ij}$ is generated by the determinantal relations 
\begin{equation}\label{quartrel}\left| \begin{array}{llll} q_{im} & q_{in} &
q_{ir} & q_{is} \\ q_{jm} & q_{jn} &
q_{jr} & q_{js}  \\ q_{km} & q_{kn} &
q_{kr} & q_{ks} \\ q_{lm} & q_{ln} &
q_{lr} & q_{ls} \end{array}\right| . \end{equation} 
\end{thm}

This theorem also holds if we have odd as well as even variables. If $U_n$ is another copy of the adjoint representation with basis $\{ u^h_n,u^x_n,u^y_n\}$ for $n\geq 0$, the supercommutative ring $$\big(Sym (\bigoplus_{n=0}^{\infty}V_n) \bigotimes \bigwedge (\bigoplus_{n=0}^{\infty}U_n) \big)^{\gs\gl_2}$$ is also generated by cubic and quadratics, as above. The only subtlety is that in the \lq\lq fermionic" determinants $c_{klm}$, the indices $k,l,m$ need not be distinct if they correspond to odd variables. For example, the following elements are nonzero: $$\left| \begin{array}{lll}a^h_k & a^x_k & a^y_k \\ u^h_l & u^x_l & u^y_l \\ u^h_l & u^x_l & u^y_l \end{array}\right|  =  2a^h_k u^x_l u^y_l + 2a^x_k u^y_l u^h_l -2 a^y_k u^x_l u^h_l, \ \ \ \ \ \ \ \ \left| \begin{array}{lll} u^h_k & u^x_k & u^y_k \\ u^h_k & u^x_k & u^y_k \\ u^h_k & u^x_k & u^y_k \end{array}\right|  = 6 u^h_k u^x_k u^y_k.$$

For the reader's convenience, we write down the vertex operators $\Theta^{\xi}_{\cW}$ and $J$ given by (\ref{defoftheta})-(\ref{defofdiff}) in terms of the basis $x,y,h$ for $\gs\gl_2$. We have
$$\Theta^x_{\cW} = 2 :\beta^x \gamma^{h'}: - :\beta^h \gamma^{y'}:    - 2:b^x c^{h'}: + :b^h c^{y'}:,$$
$$\Theta^y_{\cW} = -2 :\beta^y \gamma^{h'}: + :\beta^h \gamma^{x'}:  + 2 :b^y c^{h'}: - :b^h c^{x'}:,$$
$$\Theta^h_{\cW} = -2 :\beta^x \gamma^{x'}:  + 2 : \beta^y \gamma^{y'}: + 2 :b^x c^{x'}: - 2 :b^y c^{y'}:,$$
$$ J =  :\beta^h \gamma^{x'} c^{y'}: - :\beta^h \gamma^{y'} c^{x'}:  + 2 :\beta^x \gamma^{h'} c^{x'}: - 2 :\beta^x \gamma^{x'} c^{h'}:  - 2 :\beta^y \gamma^{h'} c^{y'}:  + 2 :\beta^y \gamma^{y'} c^{h'}:$$
$$ -:b^h c^{x'} c^{y'}: +  2 :b^x c^{x'} c^{h'}: - 2 :b^y c^{y'} c^{h'}:.$$

A natural place to look for elements of $\cW_{bas}$ is to write down the invariant normally ordered polynomials in the generators $\beta^{\xi},\gamma^{\xi'},b^{\xi}\in \cW_{hor}$, using (\ref{quadgen}), and taking into account the $\gs\gl_2$-module isomorphism $\gs\gl_2 \cong \gs\gl_2^*$. There are five quadratic, and four cubic vertex operators, given as follows:
\begin{equation}v^h =  :\beta^h \gamma^{h'}: + :\beta^x \gamma^{x'}: +:\beta^y\gamma^{y'}: ,\end{equation}
\begin{equation}v^x = \frac{1}{2}\big(:\gamma^{h'} \gamma^{h'}: + :\gamma^{x'} \gamma^{y'}:\big),\end{equation}
\begin{equation}v^y = -\frac{1}{2}\big(:\beta^h \beta^h:+ 4:\beta^x \beta^{y}:\big),\end{equation}
\begin{equation}K  = :\gamma^{h'} b^h: + :\gamma^{x'} b^x: + :\gamma^{y'} b^y:,\end{equation}
\begin{equation}Q^{\beta b} = :\beta^h b^h:+ 2 :\beta^x b^y:+2:\beta^y b^x:,\end{equation} 
\begin{equation}C^{\beta\gamma b} = -:\beta^h\gamma^{x'} b^x: + :\beta^h \gamma^{y'} b^y: - 2 :\beta^x\gamma^{h'} b^y: + :\beta^x \gamma^{x'} b^h: + 2:\beta^y \gamma^{h'} b^x: - :\beta^y \gamma^{y'} b^h:, \end{equation} 
\begin{equation}C^{\gamma bb} = - : \gamma^{h'} b^x b^y: + \frac{1}{2} :\gamma^{x'} b^x b^h:  - \frac{1}{2} :\gamma^{y'} b^y b^h:,\end{equation}
\begin{equation}C^{\beta bb} = :\beta^h b^x b^y: + :\beta^x b^y b^h: - :\beta^y b^x b^h:,\end{equation}
\begin{equation}C^{bbb} = :b^x b^y b^h:.\end{equation}
Some of these vertex operators (and similar $\gs\gl_2$-invariant vertex operators involving $c^{x'}$, $c^{y'}$, and $c^{h'}$) were written down by Akman in \cite{A}, in her study of the semi-infinite cohomology $H^{\infty + *}(\hat{\gs\gl}_2,\cS) = H^*(\cW,J(0))$. Note that $C^{\gamma bb}$ coincides with the element $C$ appearing in (\ref{defofvirasoro}) and $C^{\beta\gamma b}$ coincides with the element $H$ given by (\ref{defofh}), in the case $G=SU(2)$. Moreover, $v^x$, $v^y$, and $v^h$ coincide with the vertex operators given by (\ref{basicslii}) in the case $G=SU(2)$.

\begin{thm} \label{universalcurrent} The elements $v^x$, $v^y$, and $v^h$ are strong generators for the subalgebra $\cS^{\gs\gl_2[t]}\subset \cW_{bas}$. \end{thm}

\begin{proof} We have a ring isomorphism $$gr(\cS)^{\gs\gl_2[t]} = \cO(J_{\infty}(\gs\gl_2\oplus \gs\gl_2))^{\gs\gl_2[t]},$$ since the adjoint and coadjoint representations of $\gs\gl_2$ are isomorphic. Theorem \ref{jet} implies that $gr(\cS)^{\gs\gl_2[t]}$ is generated as a differential algebra by $\cO(\gs\gl_2\oplus \gs\gl_2)^{\gs\gl_2}$, which is a polynomial algebra on three generators. In terms of the generators $\{\beta^{\xi}_k, \gamma^{\xi'}_k|\xi = x,y,h,\ k\geq 0\}$ of $gr(\cS)$, the generators of $gr(\cS)^{\gs\gl_2[t]}$ as a differential algebra are $$\beta^h_0 \gamma^{h'}_0 + \beta^x_0 \gamma^{x'}_0 +\beta^y_0\gamma^{y'}_0,\ \ \ \ \ \ \ \ \ \ \frac{1}{2}(\gamma^{h'}_0 \gamma^{h'}_0 + \gamma^{x'}_0 \gamma^{y'}_0),\ \ \ \ \ \ \ \ \ \ -\frac{1}{2}(\beta^h_0 \beta^{h}_0 + 4\beta^x_0 \beta^{y}_0).$$  But these correspond precisely to the vertex operators $v^h$, $v^x$, and $v^y$ under the projection $\phi_2: \cS_{(2)} \rightarrow \cS_{(2)}/ \cS_{(1)} \subset gr(\cS)$. This shows that the map (\ref{defofgamma}) is surjective, and hence is an isomorphism. \end{proof}

\begin{remark} \label{uca} Since $O(\gs\gl_2,-\frac{3}{8}\kappa)$ is a simple vertex algebra, it follows that $\cS^{\gs\gl_2[t]}\cong O(\gs\gl_2,-\frac{3}{8}\kappa)$. In particular, there are no nontrivial normally ordered polynomial relations among $v^x,v^y,v^h$ and their derivatives. \end{remark}

\begin{remark} In \cite{LL}, the copy of $O(\gs\gl_2,-\frac{3}{8}\kappa)$ generated by $v^x, v^y, v^h$ was denoted by $\cA$. The main result of \cite{LL} was that $Com(\cA,\cS) = \Theta_{\cS}$, where $\Theta_{\cS}$ is the copy of $O(\gs\gl_2,-\kappa)$ generated by $\Theta^x_{\cS}, \Theta^y_{\cS}, \Theta^h_{\cS}$, which are given by (\ref{defoftheta}). Theorem \ref{universalcurrent} shows that $Com(\Theta_{\cS},\cS) = \cA$, so $\Theta_{\cS}$ and $\cA$ form a Howe pair (i.e., a pair of mutual commutants) inside $\cS$. \end{remark}

Next, an OPE calculation shows that $K$, $Q^{\beta b}$, $C^{\gamma bb}$, $C^{\beta bb}$, and $C^{bbb}$ all lie in $\cW_{bas}$. However, $C^{\beta\gamma b}$ does not lie in $\cW_{bas}$ since $$\Theta^{\xi}_{\cW}(z) C^{\beta\gamma b}(w)\sim 4b^{\xi}(w)(z-w)^{-2},\ \ \ \ \ \xi = x,y,h.$$ The image of $C^{\beta\gamma b}$ in $gr(\cW_{hor})$ lies in $(gr(\cW_{hor}))^{\gs\gl_2[t]}$, but there is no \lq\lq quantum correction" $\omega$ of lower polynomial degree, such that $C^{\beta\gamma b} +\omega \in \cW_{bas}$. This shows that the map $ gr(\cW_{bas}) \rightarrow (gr(\cW_{hor}))^{\gs\gl_2[t]}$ fails to be surjective, so we cannot reconstruct $\cW_{bas}$ from $(gr(\cW_{hor}))^{\gs\gl_2[t]}$ in a naive way, using Lemma \ref{reconlem}.

Define $\cC$ to be the vertex subalgebra of $\cW_{bas}$ generated by $v^x$, $v^y$, $v^h$, $K$, $Q^{\beta b}$, $C^{\gamma bb}$, $C^{\beta bb}$, and $C^{bbb}$. Clearly $C^{bbb}$ commutes with the other generators of $\cC$, and the following OPE relations are easy to verify:
\begin{equation}\label{kqope} K(z) Q^{\beta b} (w) \sim 0,\end{equation}
\begin{equation}v^h(z) K(w) \sim K(w)(z-w)^{-1},~~~~~v^x(z) K(w) \sim 0,~~~~~ v^y(z) K(w) \sim -Q^{\beta b} (w)(z-w)^{-1},\end{equation}
\begin{equation} \label{vcope}v^h(z) Q^{\beta b} (w) \sim -Q^{\beta b} (w)(z-w)^{-1},~~~~v^x(z) Q^{\beta b} (w) \sim -K(w)(z-w)^{-1},~~~~~ v^y(z) K(w) \sim 0.\end{equation}
\begin{equation}v^h(z) C^{\gamma bb}(w) \sim C^{\gamma bb}(w)(z-w)^{-1},~~~~v^x(z) C^{\gamma bb}(w) \sim 0,~~~~~ v^y(z) C^{\gamma bb}(w) \sim C^{\beta bb}(w)(z-w)^{-1},\end{equation}
\begin{equation}v^h(z) C^{\beta bb}(w) \sim -C^{\beta bb}(w)(z-w)^{-1},~~~~v^x(z) C^{\beta bb}(w) \sim  C^{\gamma bb}(w)(z-w)^{-1},~~~~~ v^y(z) C^{\beta bb}(w) \sim 0,\end{equation}
\begin{equation}K(z)C^{\beta bb}(w) \sim -3 C^{bbb}(w)(z-w)^{-1},~~~~~ Q^{\beta b}(z) C^{\gamma bb}(w) \sim -3 C^{bbb}(w)(z-w)^{-1},\end{equation}
\begin{equation}\label{ccope}C^{\gamma bb}(z) C^{\beta bb}(w) \sim 0,\ \ \ \ C^{\gamma bb}(z) C^{\gamma bb}(w) \sim 0, \ \ \ \ C^{\beta bb}(z) C^{\beta bb}(w) \sim 0.\end{equation}
Note that each term in the OPE of any pair of these generators is linear in the generators, so $\cC$ is a homormophic image of a current algebra associated to an $8$-dimensional Lie superalgebra $\gs$. We can obtain $\gs$ as an extension of $\gs\gl_2$ in two steps. First, let $M$ denote the irreducible, $2$-dimensional $\gs\gl_2$-module with basis $m_{-1},m_1$, regarded as an odd vector space. Define a Lie superalgebra $\tilde{\gs} = \tilde{\gs}_0 \oplus \tilde{\gs}_1$ to be the semidirect product $\gs\gl_2 \triangleright M$. In other words, $\tilde{\gs}_0 = \gs\gl_2$ and $\tilde{\gs}_{1} = M$, and the following super-commutators are satisfied: $$[h,m_{-1}] = - m_{-1},\ \ \ [h,m_1] = m_1,\ \ \ [x,m_{-1}] = m_1, \ \ \ [x,m_1] = 0, \ \ \ [y,m_{-1}] = 0,\ \ \ [y,m_1] = m_{-1},$$ $$[m_1,m_1]=0,\ \ \ \ \ [m_1,m_{-1}] = 0, \ \ \ \ \ [m_{-1},m_{-1}] = 0.$$ Next, let $N = N_0 \oplus N_1$ be a super $\tilde{\gs}$-module, where $N_0$ is a copy of the $2$-dimensional $\gs\gl_2$-module with basis $n_{-1}, n_{1}$ , $N_1$ is a copy of the trivial, 1-dimensional $\gs\gl_2$-module with basis $n_0$. The action of $\tilde{\gs}_0$ on $N$ is given as follows:
$$[h,n_{-1}] = - n_{-1},\ \ \ [h,n_1] = n_1,\ \ \ [x,n_{-1}] =  n_1, \ \ \ [x,n_1] = 0, \ \ \ [y,n_{-1}] = 0,\ \ \ [y,n_1] = n_{-1},$$ $$ [h,n_0] = 0,\ \ \ \ \ [x,n_0] = 0, \ \ \ \ \ [y,n_0] = 0,$$ $$ [m_1,n_1]=0, \ \ \  [m_1,n_{-1}] = -3 n_0,  \ \ \ [m_{-1}, n_1] = -3 n_0, \ \ \   [m_1, n_0] = 0,   \ \ \ [m_{-1}, n_0] = 0.$$
Define $\gs$ to be the semidirect product Lie superalgebra $\tilde{\gs} \triangleright N$, and define a symmetric, invariant bilinear form $B$ on $\gs$ by extending the form $-\frac{3}{8}\kappa$ on the subalgebra $\tilde{\gs}_0\subset \gs$ trivially. Let $\cF$ denote the corresponding super current algebra $O(\gs,B)$, with generators $X^u$, $u\in \gs$. By (\ref{kqope})-(\ref{ccope}) and the fact that $v^x, v^y, v^h$ generate a copy of $O(\gs\gl_2,-\frac{3}{8}\kappa)$, the map $\cF\rightarrow \cC$ sending \begin{equation}\label{supercurrenti} X^h\mapsto v^h, \ \ \ \ \ X^x\mapsto v^x,\ \ \ \ \ X^y\mapsto v^y,\ \ \ \ \ X^{m_1}\mapsto K, \ \ \ \ \ X^{m_{-1}} \mapsto Q^{\beta b}, \end{equation} \begin{equation}\label{supercurrentii}X^{n_1}\mapsto C^{\gamma bb}, \ \ \ \ \ X^{n_{-1}}\mapsto C^{\beta bb}, \ \ \ \ \ X^{n_0}\mapsto C^{bbb},\end{equation} is a vertex algebra homomorphism. 

Using the classical relation (\ref{firstrel}), we can write down certain normally ordered polynomial relations among the generators of $\cC$:
\begin{equation}  \label{ideali} :Q^{\beta b} C^{\gamma bb}:  +  :v^h C^{bbb}: +  \partial C^{bbb},\end{equation}
\begin{equation}  \label{idealii}:K C^{\gamma bb}: + 2 :v^x  C^{bbb}:,\end{equation}
\begin{equation} \label{idealiii} :Q^{\beta b} C^{\beta bb}: + 2 : v^y C^{bbb}:, \end{equation}
\begin{equation} \label{idealiv} :K C^{\beta bb}: - :v^h C^{bbb}: + \partial C^{bbb}.\end{equation} On the other hand, it is immediate from Theorem \ref{weyl} that there are {\it no} nontrivial relations in $\cW$ of the form \begin{equation}\label{forbiddeni}:f C^{\gamma bb}: + :g C^{\beta bb}: + \cdots, \end{equation} where $f,g$ are homogeneous, normally ordered polynomials in $v^x,v^y,v^h$ of the same degree, and $(\cdots)$ denotes terms in $\cW$ of lower polynomial degree. Note, however, that there does exist a relation $$:v^y C^{\gamma bb}: - \frac{1}{2}: v^h C^{\beta bb}:  + \frac{1}{4} :Q^{\beta b} C^{\beta\gamma b}: - \frac{1}{2}\partial C^{\beta bb} - \tilde{C},$$ where $\tilde{C} =  :(\partial \beta^h) b^x b^y:+ :(\partial \beta^x) b^y b^h: - :(\partial \beta^y) b^x b^h:$. This does not contradict the above statement, however, because the terms $:Q^{\beta b} C^{\beta\gamma b}:$, $v^y C^{\gamma bb}:$ and $: v^h C^{\beta bb}:$ all have polynomial degree $5$. Similarly, it follows from (\ref{quartrel}) that there are no nontrivial relations in $\cW$ of the form
\begin{equation}\label{forbiddenii}:f K: + :g Q^{\beta b}: + \cdots,\end{equation} where $f,g$ are homogeneous, normally ordered polynomials in $v^x,v^y,v^h$ of the same degree, and $(\cdots)$ denotes terms of lower polynomial degree. These observations will be important later.

Next, we observe that $\cC$ is actually a subcomplex of $\cW_{bas}$ under the differential $D(0)$. Since $K$ lies in $\cC$, it is enough to show that $J(0)$ preserves $\cC$, which is clear from the following OPEs:
\begin{equation} \label{opejkq} J(z) K(w) \sim 0,\ \ \ \ \ J(z) Q^{\beta b}(w) \sim 0, \ \ \ \ \ J(z) v^u(w) \sim 0,\ \ \ \ u = x,y,h,\end{equation}
\begin{equation} \label{opejci}J(z) C^{\gamma bb}(w) \sim K(w)(z-w)^{-2} + (: v^h K: -2 :v^x Q^{\beta b}:~ -\partial K)(z-w)^{-1},\end{equation}
\begin{equation} \label{opejcii} J(z) C^{\beta bb}(w) \sim -Q^{\beta b} (w)(z-w)^{-2} + (2: v^y K:~+~ :v^h Q^{\beta b} :~ + \partial Q^{\beta b} )(z-w)^{-1},\end{equation}
\begin{equation} \label{opejciii} J(z) C^{bbb}(w) \sim -~:K Q^{\beta b} :(z-w)^{-1}.\end{equation} It is amusing to note that $J(0)$ acts {\it nonlinearly} on the generators of $\cC$, so its action on $\cC$ does not come from a differential Lie superalgebra structure on $\gs$. We conjecture that $\cW_{bas}=\cC$, although we are unable to prove this at present. Let $\cC^{(k)} = \cC \cap \cW_{bas}^{(k)}$, which is the homogeneous subspace of $\cC$ of $b$-number $k$. Since $\cW_{bas}^{(0)} = \cS^{\gs\gl_2[t]} = \cC^{(0)}$ by Theorem \ref{universalcurrent}, our conjecture clearly holds for $b$-number zero. It will suffice for our purposes to show that $\cW_{bas}^{(1)} = \cC^{(1)}$ as well. First, we need a characterization of $\cW_{bas}^{(1)}$ which is analogous to the characterization of $\cS^{\gs\gl_2[t]}$ given by Lemma \ref{bnumbzero}.

\begin{lemma} \label{bnumbone} We have $\cW_{bas}^{(1)}= \cW_{hor}^{(1)}\cap Ker(J(0))$. \end{lemma}

\begin{proof} For any simple Lie algebra $\gg$, the inclusion \begin{equation} \label{charwbi} \cW(\gg)_{bas}^{(1)} \subset \cW(\gg)_{hor}^{(1)}\cap Ker(J(0))\end{equation} is equivalent to the injectivity of the map \begin{equation}\label{classinj} \cS(\gg)^{\gg[t]}\rightarrow H^*(\cW(\gg)_{bas},J(0))\end{equation} given by Lemma \ref{bnumbzero}. In the case $\gg = \gs\gl_2$, we have $\cS^{\gs\gl_2[t]} \cong O(\gs\gl_2,-\frac{3}{8}\kappa)$, which is a simple vertex algebra, so the injectivity of $\cS^{\gs\gl_2[t]}\rightarrow H^*(\cW_{bas},J(0))$ is apparent. We expect (although we are unable to prove) that the map (\ref{classinj}) is injective for any simple $\gg$. 

Next, we claim that the opposite inclusion $Ker(J(0)) \cap \cW(\gg)_{hor}^{(1)} \subset \cW(\gg)_{bas}$ holds for any simple $\gg$. First, any $\omega\in \cW(\gg)_{hor}^{(1)}$ can be written in the form
$$\omega = \sum_{k=0}^l \sum_{j=1}^n :\partial^k b^{\xi_j} P_k^{\xi_j}: ,$$ for some fixed $l$. In this notation, $\{\xi_1,\dots,\xi_n\}$ is an orthonormal basis for $\gg$ relative to the Killing form, and each $P^{\xi_j}_k$ lies in $\cS(\gg)$.

We need the following computations: $$ \big(\sum_{i=1}^n :c^{\xi_i}\Theta^{\xi_i}_{\cE}: \big)\circ_0  b^{\xi_j} =2 \Theta^{\xi_j}_{\cE},\ \ \ \ \ \ \ \frac{1}{m!} \Theta^{\xi_i}_{\cE}\circ_m \partial^l b^{\xi_j} = \binom{k}{m} \partial^{k-m}b^{[\xi_i,\xi_j]}.$$
Using this calculation, we obtain
$$J\circ_0\omega = \sum_{j=1}^n \sum_{k=0}^l :(\partial^k \Theta^{\xi_j}_{\cW})P^j_k :+ \sum_{i,j=1}^n \sum_{k=0}^l \sum_{m\geq 0} \frac{1}{m!}: (\partial^m c^{\xi_i}) (\partial^k b^{\xi_j})( \Theta^{\xi_i}_{\cS} \circ_m P^{\xi_j}_k) : .$$ The term $ \sum_{j=1}^n \sum_{k=0}^l :(\partial^k \Theta^{\xi_j}_{\cW})P^j_k :$ can be rewritten in the form $$  \sum_{j=1}^n \sum_{k=0}^l :(\partial^k\Theta^{\xi_j}_{\cS}) P^{\xi_j}_k:  + \sum_{i,j=1}^n \sum_{k=0}^l \sum_{m=0}^k \binom{k}{m}  : (\partial^{k-m} b^{[\xi_j,\xi_i]}) (\partial^{m} c^{\xi'_i})P^{\xi_j}_k: $$
$$=  \sum_{j=1}^n \sum_{k=0}^l :(\partial^k\Theta^{\xi_j}_{\cS} )P^{\xi_j}_k: + \sum_{i,j=1}^n \sum_{k=0}^l \sum_{m=0}^k \binom{k}{m}  :  (\partial^{m} c^{\xi'_i}) (\partial^{l-m} b^{[\xi_i,\xi_j]})P^{\xi_j}_k : $$ 
$$=   \sum_{j=1}^n \sum_{k=0}^l : (\partial^k\Theta^{\xi_j}_{\cS}) P^{\xi_j}_k: + \sum_{i,j=1}^n \sum_{k=0}^l \sum_{m=0}^k \frac{1}{m!} :  (\partial^{m} c^{\xi'_i}) (\Theta^{\xi_i}_{\cE} \circ_m \partial^l b^{\xi_j})P^{\xi_j}_k:  .$$ Collecting terms, we have
\begin{equation}\label{jcircow} J\circ_0 \omega=\sum_{j=1}^n \sum_{k=0}^l :(\partial^k\Theta^{\xi_j}_{\cS}) P^{\xi_j}_k: +  \frac{1}{m!} \sum_{m\geq 0} :(\partial^m c^{\xi'_i})(\Theta^{\xi_i}_{\cW}\circ_m \omega):.\end{equation} Since $J\circ_0 \omega = 0$, and the term $\sum_{j=1}^n \sum_{k=0}^l :(\partial^k\Theta^{\xi_j}_{\cS}) P^{\xi_j}_k: $ does not depend on $c^{\xi'_i}$ and its derivatives, we conclude that $\Theta^{\xi_i}_{\cW} \circ_m \omega = 0$ for all $m\geq 0$. Hence $\omega$ lies in $\cW(\gg)_{bas}$. \end{proof}

\begin{remark} The characterization $\cW_{bas}^{(k)}= \cW_{hor}^{(k)}\cap Ker(J(0))$ fails for $k\geq 2$; see for example (\ref{opejci})-(\ref{opejciii}). \end{remark}

\begin{thm} \label{bnumbonedesc} For $\gg = \gs\gl_2$, $\cW_{bas}^{(1)} = \cC^{(1)}$, and hence has a basis consisting of normally ordered polynomials in $v^x, v^h, v^y, Q^{\beta b},K$ and their derivatives, which are linear in $Q^{\beta b}$ and $K$, and their derivatives. In particular, $\cW_{bas}^{(1)}$ is homogeneous of even polynomial degree. \end{thm}

\begin{proof} Recall that $$gr(\cW_{hor}) \cong Sym  (\bigoplus_{k\geq 0}(V_k\oplus V^*_k)) \bigotimes \bigwedge(\bigoplus_{k\geq 0 } U_k),$$ where $V_k$ and $U_k$ are the copies of $\gg = \gs\gl_2$ with bases $\{\beta^{\xi}_k|~\xi=x,y,h\}$ and $\{b^{\xi}_k| \xi=x,y,h\}$, respectively, and $V^*_k$ is the copy of $\gg^*$ with basis $\{\gamma^{\xi'}_k| \xi=x,y,h\}$. We are interested in the subspace $gr(\cW_{hor}^{(1)})$ of $b$-number $1$, which is linear in the variables $b^{\xi}_k$. Since this subspace is linear in these variables, it does not matter whether they are regarded as even or odd variables. So $gr(\cW_{hor}^{(1)})$ is isomorphic to the subspace of $$Sym(\bigoplus_{k\geq 0}(V_k\oplus V^*_k\oplus U_k)) = \cO(J_{\infty}(V^*\oplus V \oplus U^*))$$ which is linear in the generators of $\cO(J_{\infty}(U^*))$. Here $V = \gg$ and $U^*,V^*$ are copies of $\gg^*$. By Theorem \ref{jet}, $\cO(J_{\infty}(V^*\oplus V \oplus U^*))^{\gs\gl_2[t]}$ is generated as a differential algebra by $\cO(V^*\oplus V \oplus U^*)^{\gs\gl_2}$, which by Theorem \ref{weyl} is generated by six quadratics
$$q_{VV},  \ \ \ \ \ \ q_{VV^*},\ \ \ \ \ \ q_{V^*V^*},\ \ \ \ \ \ q_{VU^*},\ \ \ \ \ \ q_{V^*U^*},\ \ \ \ \ \ q_{U^*U^*},$$ and one cubic $c_{VV^*U^*}$. The subspace which is linear in the generators of $\cO(J_{\infty}(U^*))$ is the $\langle \mathbb{C}[q_{VV}, q_{VV^*}, q_{V^*V^*}]\rangle $-module generated by $q_{VU^*}$, $q_{V^*U^*}$ and $c_{VV^*U^*}$, and their derivatives. (In this notation, $\langle \mathbb{C}[q_{VV}, q_{VV^*}, q_{V^*V^*}]\rangle$ denotes the algebra generated by $q_{VV}, q_{VV^*}, q_{V^*V^*}$, and their derivatives).

Recall the projections $\phi_d: \cW_{(d)} \rightarrow \cW_{(d)}/ \cW_{(d-1)} \subset gr(\cW)$. Define \begin{equation} \label{mgeni} v^h_m = \phi_2 (\partial^m v^h),\ \ \ \ \ \ \ v^x_m = \phi_2 (\partial^m v^x),\ \ \ \ \ \ \ v^y_m = \phi_2 (\partial^m v^y), \end{equation} \begin{equation} \label{mgenii}Q^{\beta b}_m= \phi_2 (\partial^m Q^{\beta b}), \ \ \ \ \ \ \  K_m= \phi_2 (\partial^m K),\ \ \ \ \ \ \ C^{\beta \gamma b}_m = \phi_3(\partial^m C^{\beta\gamma b}).\end{equation} 
Clearly $q_{VV}$, $ q_{VV^*}$, $q_{V^*V^*}$,  $q_{VU^*}$, $q_{V^*U^*}$, and $c_{VV^*U^*}$ correspond to $v^x_0$, $v^h_0$, $v^y_0$, $K_0$, $Q^{\beta b}_0$, and $C^{\beta \gamma b}_0$, respectively. It follows that $(gr(\cW_{hor}^{(1)}))^{\gs\gl_2[t]}$ is a module over $gr(\cS^{\gs\gl_2[t]}) = \mathbb{C}[v^x_m, v^y_m, v^h_m|m\geq 0]$, with generators  $K_n$, $Q^{\beta b}_n$, and $C^{\beta \gamma b}_n$, for $n\geq 0$.

Since $\cW_{bas} = ( \cW_{hor})^{\gs\gl_2[t]}$, we have an injective linear map $$gr(\cW_{bas}^{(1)}) \hookrightarrow (gr(\cW_{hor}^{(1)}))^{\gs\gl_2[t]}.$$ Recall that this map is not surjective, since $C^{\beta\gamma b}_0$  does not lie in the image, so we cannot reconstruct $\cW_{bas}^{(1)}$ in a naive way. However, for any $\omega\in \cW_{bas}^{(1)}$ of polynomial degree $d$, $\phi_d(\omega)\in gr(\cW)$ is a polynomial in $v^x_m$, $v^h_m$, $v^y_m$, $K_m$, $Q^{\beta b}_m$, and $C^{\beta \gamma b}_m$ for $m\geq 0$. It follows that $\omega$ can be written in the form $$\omega = \omega_{d} + \omega_{<},$$ where the leading term $\omega_{d}$ is a normally ordered polynomial in $v^x$, $v^y$, $v^h$, $K$, $Q^{\beta b}$, $C^{\beta\gamma b}$, and their derivatives, and $\omega_{<}\in \cW_{(d-2)}$. Since the operators $\Theta_{\cW}^{\xi}(k)$, $k\geq 0$ on $\cW_{hor}$ either preserve polynomial degree, or lower it by $2$, we have $\cW_{bas} = \cW_{bas}^{even} \oplus \cW_{bas}^{odd}$, where $\cW_{bas}^{even}$ and $\cW_{bas}^{odd}$ are homogeneous of even and odd polynomial degrees, respectively. We can therefore deal with the odd and even cases separately. Suppose first that $\omega\in \cW_{bas}^{(1)}$ and has even polynomial degree $d$. Then $\omega_{d}$ is a normally ordered polynomial in $v^x$, $v^y$, $v^h$, $K$, $Q^{\beta b}$ and their derivatives, so $\omega_{d}\in \cD$ by definition. In particular, $\omega_{d} \in \cW^{(1)}_{bas}$, so $\omega_{<}$ must also lie in $\cW^{(1)}_{bas}$. Since $\omega_{<}$ lies in $\cW_{(d-2)}$, it follows by induction on $d$ that $\omega\in \cD$.

Next, suppose that $\omega$ has odd polynomial degree $d$. Then the leading term $\omega_{d}$ is of the form $$\sum_{k=0}^l :P_k \partial^k C^{\beta\gamma b}:$$ for some fixed $l$, where each $P_k$ lies in $\cS^{\gs\gl_2[t]}$, and $P_l\neq 0$. Unlike the case where $d$ is even, $\omega_{d}$ does not lie in $\cW_{bas}^{(1)}$, since $\Theta^{\xi}_{\cW}(l+1)( \omega_{d}) = 4 (l+1)! :P_l b^{\xi}:$, for $\xi = x,y,h$. 

Recall the vertex operator $F = -(:b^x c^{x'}: + :b^y c^{y'}: +:b^h c^{h'}:)$, which is part of the TVA structure given by (\ref{newtopii}). We calculate $$F(z) C^{\beta\gamma b}(w) \sim -C^{\beta\gamma b}(w)(z-w)^{-1},$$ which implies that $F(l)(\partial^l C^{\beta\gamma b}) = -l! C^{\beta\gamma b}$, and $F(l)(\partial^k C^{\beta\gamma b}) = 0$ for $0\leq k<l$. Since $F(l)$ acts trivially on each $P_k$,  we have $$-\frac{1}{l!}F(l)(\sum_{k=0}^l :P_k \partial^k C^{\beta\gamma b}:) = :P_l C^{\beta\gamma b}:.$$ Since $F(k)(\cW_{bas})\subset \cW_{bas}$ for all $k\geq 0$, we can assume without loss of generality that $l=0$, so that $\omega_{d}$ is of the form $:P C^{\beta\gamma b}:$, with $P\in \cS^{\gs\gl_2[t]}$. 

Next, we compute $$J(0)(C^{\beta\gamma b}) = -2:v^h v^h: -8 :v^x v^y: + 4 \omega_{\cW} +2 \partial v^h,$$ where $\omega_{\cW}$ is the Virasoro element in $\cW_{(2)}$ given by (\ref{weilvir}). Since $J(0)$ annihilates $P$, we have $$J(0)(\omega_{d}) = :P(-2:v^h v^h: -8 :v^x v^y:): + \cdots,$$ where $(\cdots)$ lies in $\cW_{(d-1)}$. Hence $\phi_{d+1}(J(0)(\omega_{d}))$ is a nonzero element of $gr(\cW)$ of polynomial degree $d+1$. Since $J(0)$ can raise the polynomial degree by at most $1$, and $\omega_{<}$ has polynomial degree at most $d-2$, $J(0)(\omega_{<})$ can have polynomial degree at most $d-1$. Hence the term $:P(-2:v^h v^h: -8 :v^x v^y:):$ cannot be canceled. In particular, $J(0)(\omega)\neq 0$, which contradicts $\omega\in \cW_{bas}^{(1)}$, by Lemma \ref{bnumbone}. \end{proof}

\section{The structure of $H^*(\cW_{bas},K(0))$}

Recall that $\cW_{bas}$ is a double complex under the commuting differentials $K(0)$ and $J(0)$, whose sum is $D(0)$. Using the description of $\cW_{bas}^{(0)}$ and $\cW_{bas}^{(1)}$ given by Theorems \ref{universalcurrent} and \ref{bnumbonedesc}, we will construct an interesting subalgebra of $H^*(\cW_{bas},K(0))$. Later, we will use the results in this section together with the spectral sequence of the double complex to study $\H^*_{SU(2)}(\mathbb{C})$. First, since $K(0)(\cW_{bas}^{(k)})\subset \cW_{bas}^{(k+1)}$, we may regrade the complex $(\cW_{bas}^*,K(0))$ by $b$-number. We denote this new complex by $(\cW_{bas}^{(*)},K(0))$, and we denote its cohomology by $H^{(*)}(\cW_{bas},K(0))$. Clearly we have a linear isomorphism \begin{equation}\label{regradc} \bigoplus_{k\geq 0} H^{(k)}(\cW_{bas},K(0)) \cong \bigoplus_{l\in \mathbb{Z}} H^l (\cW_{bas},K(0)).\end{equation} 

Let $\cD$ denote the subalgebra of $\cW_{bas}$ generated by $v^x$, $v^y$, $v^h$, $K$, and $Q^{\beta b}$, which contains both $\cW_{bas}^{(0)}$ and $\cW_{bas}^{(1)}$. Clearly $(\cD^*,K(0))$ and $(\cD^{(*)},K(0))$ are subcomplexes of $(\cW_{bas}^*,K(0))$ and $(\cW_{bas}^{(*)},K(0))$, respectively, where $\cD^{(k)} = \cD\cap \cW_{bas}^{(k)}$. There is a similar linear isomorphism $$ \bigoplus_{k\geq 0} H^{(k)}(\cD,K(0)) \cong \bigoplus_{l\in \mathbb{Z}} H^l (\cD,K(0)).$$

Let $\langle\gamma\rangle$ denote the (abelian) subalgebra of $\cW$ generated by $\gamma^{x'}, \gamma^{y'},\gamma^{h'}$. It is easy to see from Theorem \ref{universalcurrent} that $\langle \gamma\rangle^{\gs\gl_2[t]}$ is just the polynomial algebra generated by $\partial^k v^x$ for $k\geq 0$.

\begin{lemma} \label{dcoh} $H^{(k)}(\cD,K(0))$ vanishes for $k>0$, and $H^{(0)}(\cD,K(0))$ is isomorphic to $\langle \gamma\rangle^{\gs\gl_2[t]}$. \end{lemma}

\begin{proof} Since $K(0)$ preserves the filtration (\ref{filtw}), $K(0)$ acts on $gr(\cD)$, and we may first consider the cohomology $H^{(k)}(gr(\cD),K(0))$. Recall from (\ref{mgeni})-(\ref{mgenii}) that $v^x_m$, $v^y_m$, $v^h_m$, $Q^{\beta b}_m$, $K_m$ denote the images of $\partial^m v^x$, $\partial^m v^y$, $\partial^m v^h$, $\partial^m Q^{\beta b}$, $\partial^m K$ in $gr(\cD)$ under the map $\phi_2: \cD_{(2)} \rightarrow \cD_{(2)} / \cD_{(1)} \subset gr(\cD)$. Clearly $\{v^x_m, v^y_m,v^h_m, Q^{\beta b}_m, K_m|~ m\geq 0\}$ generates $gr(\cD)$ as a ring. We will compute $H^{(k)}(gr(\cD),K(0))$ by constructing a contracting homotopy on this complex. The action $K(0)$ on $gr(\cW_{hor})$ is defined on the generators $\{\beta^{\xi}_k, \gamma^{\xi'}_k,b^{\xi}_k|\ \xi=x,y,h,\ \ k\geq 0\}$ by $$\beta^{\xi}_k\mapsto -b^{\xi}_k, \ \ \ \ \ \gamma^{\xi'}\mapsto 0,\ \ \ \ \  b^{\xi}_k\mapsto 0.$$ From this, we see that on the generators of $gr(\cD)$, $K(0)$ acts as follows: $$K(0)(v^y_m) = Q^{\beta b}_m,\ \ \ \ K(0)(v^h_m) = -K_m, \ \ \ \ K(0)(v^x_m) = 0,\ \ \ \ K(0)(K_m) = 0,\ \ \ \ K(0)(Q^{\beta b}_m) = 0.$$ 
Define a vertex operator $R = :\beta^h c^{h'}: + :\beta^x c^{x'}: + :\beta^y c^{y'}:$, and let $\tilde{R}(0)$ denote the linear map $$\sum_{k\geq 0} \beta^h(-k-1) c^{h'}(k) + \beta^x (-k-1)c^{x'}(k) + \beta^y(-k-1) c^{y'}(k) \in End(\cW_{hor}),$$ which is the sum of all terms of the Fourier mode $R(0)$ which only contain annihilation modes of $c^{\xi'}$ for $\xi=x,y,h$. There is an induced derivation on $gr(\cW_{hor})$ which we also denote by $\tilde{R}(0)$, defined on generators by
$$b^{\xi}_k\mapsto \beta^{\xi}_k, \ \ \ \ \ \gamma^{\xi'}\mapsto 0,\ \ \ \ \  \beta^{\xi}_k\mapsto 0.$$ It follows that $\tilde{R}(0)$ acts on $gr(\cD)$ as follows:
$$\tilde{R}(0)(Q^{\beta b}_m) = -2v^y_m,\ \ \ \ \tilde{R}(0)(K_m) = v^h_m, \ \ \ \ \tilde{R}(0)(v^x_m) = 0,\ \ \ \ \tilde{R}(0)(v^h) = 0,\ \ \ \ \tilde{R}(0)(v^y) = 0.$$ Moreover, as derivations on $gr(\cW_{hor})$ we compute $[K(0),\tilde{R}(0)] = S$, where $S$ is the diagonalizable operator defined on generators by:
$$S(\beta^{\xi}_m) = -\beta^{\xi}_m,\ \ \ \ \ S(b^{\xi}_m) = -b^{\xi}_m,\ \ \ \ \ S(\gamma^{\xi'}_m) = 0.$$
It follows that $$S(v^x_m) = 0,\ \ \ \ S(v^y_m) = -2 v^y_m,\ \ \ \ S(v^h_m) =- v^h_m,\ \ \ \ S(Q^{\beta b}_m) = -2 Q^{\beta b}_m,\ \ \ \ S(K_m) = -K_m,$$ so $\tilde{R}(0)$ is a contracting homotopy for $K(0)$ outside the (trivial) subcomplex generated by $\{v^x_m|\ m\geq 0\}$. Hence $H^{(k)}(gr(\cD),K(0)) \cong \langle \gamma \rangle^{\gs\gl_2[t]}$ for $k=0$, and vanishes for $k>0$.

Finally, we need to compare $H^{(k)}(gr(\cD),K(0))$ with $H^{(k)}(\cD,K(0))$. Since $K(0)$ preserves polynomial degree, $H^{(k)}(gr(\cD),K(0))$ has a grading $$H^{(k)}(gr(\cD),K(0)) = \bigoplus_{d\geq 0} H^{(k)}(gr(\cD),K(0))^d,$$ where $H^{(k)}(gr(\cD),K(0))^d$ is homogeneous of polynomial degree $d$. Similarly, $H^{(k)}(\cD,K(0))$ admits a filtration $$H^{(k)}(\cD,K(0))_{0}\subset H^{(k)}(\cD,K(0))_{1}\subset H^{(k)}(\cD,K(0))_{2}\subset \cdots,$$ where $H^{(k)}(\cD,K(0))_{d}$ denotes the subspace of $H^{(k)}(\cD,K(0))$ admitting a representative of polynomial degree at most $d$. 
There is an injective linear map $$H^{(k)}(\cD,K(0))_{d} / H^{(k)}(\cD,K(0))_{d-1} \rightarrow H^{(k)}(gr(\cD),K(0))^d.$$ It follows that $H^{(k)}(\cD,K(0))$ vanishes for $k>0$. Since $gr(\cD^{(0)})\cong \cD^{(0)} \cong \langle \gamma \rangle^{\gs\gl_2[t]}$, the claim follows. \end{proof}

By Theorems \ref{universalcurrent} and \ref{bnumbonedesc}, both $\cW_{bas}^{(0)}$ and $\cW_{bas}^{(1)}$ lie in $\cD$. It follows from Lemma \ref{dcoh} that $H^{(0)}(\cW_{bas},K(0)) \cong \langle \gamma \rangle^{\gs\gl_2[t]}$ and $H^{(1)}(\cW_{bas},K(0)) = 0$. However, we cannot compute $H^{(k)}(\cW_{bas},K(0))$ for $k\geq 2$ using this approach, because $gr(\cW_{bas})$ is not preserved by $\tilde{R}(0)$. For example, let $C^{\gamma bb}_0$ and $C^{\beta\gamma b}_0$ denote the images of $C^{\gamma bb}$ and $C^{\beta \gamma b}$ in $gr(\cW)$. Clearly $C^{\gamma bb}_0\in gr(\cW_{bas})$, but $\tilde{R}(0)(C^{\gamma bb}_0) = C^{\beta\gamma b}_0$, which does not live in $gr(\cW_{bas})$. As we shall see, $H^{(2)}(\cW_{bas},K(0))$ has a very rich structure. 

\begin{lemma} \label{oddpolydeg} The map $H^{(2)}(\cC,K(0)) \rightarrow H^{(2)}(\cW_{bas},K(0))$ induced by the inclusion $\cC\hookrightarrow \cW_{bas}$ is injective. Moreover, any nonzero element $\omega\in \cW_{bas}^{(2)}\cap Ker(K(0))$ of odd polynomial degree represents a nontrivial class in $H^{(2)}(\cW_{bas},K(0))$. By (\ref{regradc}), $\omega$ then represents a nontrivial class in $H^*(\cW_{bas},K(0))$ as well. \end{lemma}

\begin{proof} The first statement is clear since $\cW_{bas}^{(1)} = \cC^{(1)}$. The second statement is clear because there are no elements of $\cW_{bas}^{(1)}$ of odd polynomial degree, and $K(0)$ preserves the parity of polynomial degree. \end{proof}

Hence we can construct nontrivial elements of $H^{(2)}(\cW_{bas},K(0))$ by finding elements of $\cC^{(2)}\cap Ker(K(0))$ of odd polynomial degree. Let $\cI$ denote the kernel of the map $\cF\rightarrow \cC$ given by (\ref{supercurrenti})-(\ref{supercurrentii}). Since $\cF$ is a complex under the differential $X^{m_1}(0)$ corresponding to $K(0)$, we obtain a short exact sequence of complexes\begin{equation} \label{shortexact}0\rightarrow \cI \rightarrow \cF \rightarrow \cC \rightarrow 0.\end{equation} Note that $C^{\gamma bb}$ represents a nontrivial class in $H^{(2)}(\cC,K(0))$, since it has polynomial degree $3$. In fact, $C^{\gamma bb}$ corresponds to $X^{n_1}\in \cF$, which is easily seen to represent a nontrivial class in $H^*(\cF ,X^{m_1}(0))$. So the class of $C^{\gamma bb}$ in $H^{(2)}(\cC,K(0))$ lies in the image of the map $H^*(\cF,X^{m_1}(0))\rightarrow H^*(\cC,K(0))$.

We can obtain new classes in $H^*(\cC,K(0))$ which correspond to elements of $H^*(\cI,X^{m_1}(0))$ via the connecting homomorphism in the long exact sequence of (\ref{shortexact}). In other words, we seek elements of $\cF$ which do not map to zero under $X^{m_1}(0)$, but rather, map to elements of $\cI$. For example, consider $$:X^y X^{n_1}:~ - \frac{1}{4} : X^h X^{n_{-1}}:~  - \frac{5}{12} \partial X^{n_{-1}} \in \cF.$$ Under $X^{m_1}(0)$, it maps to $$: X^{m_{-1}}  X^{n_1}:  +  :X^h X^{n_0}: +  \partial X^{n_0},$$ which is the element of $\cI$ corresponding to the relation (\ref{ideali}). It follows that the corresponding element \begin{equation}\label{defhiv} h_4 = ~:v^y C^{\gamma bb}:~ - \frac{1}{4} : v^h C^{\beta bb}:~  - \frac{5}{12} \partial C^{\beta bb}\in \cC^{(2)} \end{equation} lies in the kernel of $K(0)$. Since $h_4$ has the form (\ref{forbiddeni}), it is nonzero, and since it has polynomial degree $5$, it represents a nontrivial class in $H^{(2)}(\cW_{bas},K(0))$. Note that $h_4$ has weight $4$ and cohomology degree $-4$, so it represents a nontrivial class in $H^{-4}(\cW_{bas},K(0))$.

More generally, for $n\geq 2$ define \begin{equation} \label{hgenformula} h_{2n+2} = : (v^y)^n C^{\gamma bb}:~ -\frac{n}{2n+2} : (v^y)^{n-1} v^h C^{\beta bb}: - \frac{n^2-n}{2n^2+3n+1} :\partial v^y (v^y)^{n-2} C^{\beta bb}: \end{equation}
$$  - \frac{2n^2+3n}{4n^2+6n+2}  : (v^y)^{n-1}\partial C^{\beta bb}:.$$
\begin{lemma} \label{genkclass}For all $n\geq 2$, $h_{2n+2}$ represents a nontrivial class in $H^{-4n}(\cW_{bas},K(0))$ of conformal weight $2n+2$. \end{lemma}
\begin{proof} Clearly $h_{2n+2}$ is homogeneous of weight $2n+2$ and degree $-4n$. Using the derivation property of $K(0)$ and the relations (\ref{ideali})-(\ref{idealiv}), the following calculations show that $h_{2n+2}$ lies in the kernel of $K(0)$:

\begin{equation}K(0) \big( : (v^y)^n C^{\gamma bb}:\big) = n : (v^y)^{n-1} Q^{\beta b} C^{\gamma bb}: = -n :(v^y)^{n-1} v^h C^{bbb}: - n :(v^y)^{n-1} \partial C^{bbb}:,\end{equation}

\begin{equation}K(0) \big( : (v^y)^{n-1} v^h C^{\beta bb}:\big) = (n-1) :(v^y)^{n-2}  Q^{\beta b} v^h C^{\beta bb}: - :(v^y)^{n-1} K C^{\beta bb}: - 3:(v^y)^{n-1} v^h C^{bbb}: \end{equation} $$ = (n-1) :(v^y)^{n-2} v^h Q^{\beta b} C^{\beta bb}: + (n-1) :(v^y)^{n-2} (\partial Q^{\beta b})C^{\beta bb}:- :(v^y)^{n-1} K C^{\beta bb}: - 3:(v^y)^{n-1} v^h C^{bbb}:$$ $$= -2(n-1) :(v^y)^{n-1} v^h C^{bbb}: +4(n-1):(v^y)^{n-2} (\partial v^y) C^{bbb}: + (n-1):(v^y)^{n-2} (\partial Q^{\beta b}) C^{\beta bb}:  $$ $$- :(v^y)^{n-1} v^h C^{bbb}:  + :(v^y)^{n-1} \partial C^{bbb}:  - 3:(v^y)^{n-1} v^h C^{bbb}:$$ $$ = (-2n-2) :(v^y)^{n-1} v^h C^{bbb}: +4(n-1):(v^y)^{n-2} (\partial v^y) C^{bbb}: $$ $$+ (n-1):(v^y)^{n-2} (\partial Q^{\beta b}) C^{\beta bb}: + :(v^y)^{n-1} \partial C^{bbb}:,$$

\begin{equation}K(0) \big( : (v^y)^{n-1} \partial C^{\beta bb}:\big) = (n-1) :(v^y)^{n-2} Q^{\beta b} \partial C^{\beta bb}: - 3 :(v^y)^{n-1} \partial C^{bbb}: \end{equation}$$ = (n-1) :(v^y)^{n-2} \partial \big ( Q^{\beta b} C^{\beta bb}\big): - (n-1) :(v^y)^{n-2} (\partial Q^{\beta b}) C^{\beta bb}: - 3 :(v^y)^{n-1} \partial C^{bbb}:$$ $$ =   -2(n-1) :(v^y)^{n-2} \partial \big(v^y C^{bbb}\big): - (n-1) :(v^y)^{n-2} (\partial Q^{\beta b}) C^{\beta bb}: - 3 :(v^y)^{n-1} \partial C^{bbb}:$$ $$ = (-2n-1) :(v^y)^{n-1} \partial C^{bbb}: - 2(n-1) :(v^y)^{n-2} (\partial v^y) C^{bbb}: - (n-1) :(v^y)^{n-2} (\partial Q^{\beta b}) C^{\beta bb}:,$$

\begin{equation}K(0) \big( : (v^y)^{n-2} \partial v^y C^{\beta bb}:\big)\end{equation} $$= (n-2) :(v^y)^{n-3}  Q^{\beta b} (\partial v^y) C^{\beta bb}: + :(v^y)^{n-2} ( \partial  Q^{\beta b}) C^{\beta bb}: - 3 :(v^y)^{n-2} ( \partial v^y) C^{bbb}: $$ $$= (-2n+1) :(v^y)^{n-2} (\partial v^y) C^{bbb}: + :(v^y)^{n-2} (\partial  Q^{\beta b}) C^{\beta bb}:.$$

Finally, since $h_{2n+2}$ has odd polynomial degree $2n+3$, it represents a nontrivial class in $H^{(2)}(\cW_{bas},K(0))$. It follows that $h_{2n+2}$ represents a nontrivial class in $H^{-4n}(\cW_{bas},K(0))$ as well. \end{proof}

Next, we construct additional nontrivial classes in $H^*(\cW_{bas},K(0))$ by making use of its algebraic structure. Since $v^x$ represents a class in $H^4(\cW_{bas},K(0))$ of weight zero, we can obtain new elements by taking circle products of $v^x$ with $h_{2n+2}$. For example, define $$f_3 = v^x\circ_0 h_4= :v^h C^{\gamma bb}: + \frac{2}{3} :v^x C^{\beta bb}:-\frac{5}{3} \partial C^{\gamma bb} \in  \cW_{bas}^{(2)} \cap Ker(K(0)). $$ Since $f_3$ has polynomial degree $5$, it represents a nontrivial class in $H^0(\cW_{bas},K(0))$ of weight $3$. 

\begin{lemma} \label{monomials} For any integers $n\geq 1$ and $0\leq d\leq n$, starting from the normally ordered monomial $\mu = :(v^y)^{n-d} (v^h)^{d}:$ of polynomial degree $2n$, there exists a vertex operator of the form \begin{equation} \label{specialvo} :\mu C^{\gamma bb}: + :f C^{\beta bb}: + \cdots \ \ \in \cW_{bas}^{(2)} \cap Ker(K(0))\end{equation} of polynomial degree $2n+3$, which represents a nontrivial class in $H^{-4(n-d)}(\cW_{bas},K(0))$ of weight $2n-d+2$. Here $f$ is a normally ordered polynomial in $v^x, v^y, v^h$ of polynomial degree $2n$, and the term $(\cdots)$ has polynomial degree at most $2n+1$. \end{lemma}

\begin{proof} In fact, we will prove slightly more. Not only does the vertex operator (\ref{specialvo}) exist for all $n$ and $d$, but in the case $d<n$, we will show that there appears in $f$ a monomial of the form $:(v^y)^{n-d-1}(v^h)^{d+1}:$ whose coefficient $\lambda$ lies in the open interval $(-1,0)$. 

By Lemma \ref{genkclass}, the vertex operator (\ref{specialvo}) exists when $d=0$ for all $n\geq 1$. Moreover, the coefficient of $:(v^y)^{n-1}(v^h):$ in $f$ is $-\frac{n}{2n+2}$ by (\ref{hgenformula}), which certainly lies in $(-1,0)$. Assume inductively that the statement (together with our additional assumption on $\lambda$) holds for $n-1$ and all $d = 0,\dots,n-1$. Since it holds for $n$ at the value $d=0$, we may proceed by induction on $d$, so we assume it holds for $d-1$.

Thus our inductive assumption is that there exists a normally ordered polynomial in $\cW_{bas}^{(2)}\cap Ker(K(0))$ of the form $$:(v^y)^{n-d+1} (v^h)^{d-1} C^{\gamma bb}: + \lambda : (v^y)^{n-d}(v^h)^{d} C^{\beta bb}: + \cdots,\ \ \ \ \ \ \ \lambda\in (-1,0)$$ where $(\cdots)$ consists of terms which either have polynomial degree at most $2n+1$, or do not depend on $C^{\gamma bb}$, and the term $ : (v^y)^{n-d}(v^h)^{d} C^{\beta bb}: $ does not appear in $(\cdots)$. Apply the operator $v^x\circ_0$, which preserves $\cW_{bas}^{(2)} \cap Ker(K(0))$. Using (\ref{vcope}) and the derivation property of $v^x\circ_0$, we obtain the following expression:
\begin{equation}\label{intercrap} (n-d+1) :(v^y)^{n-d} (v^h)^{d} C^{\gamma bb}: - 2(d-1) :(v^y)^{n-d+1} (v^h)^{d-2}v^x C^{\gamma bb} :\end{equation} $$+ \lambda \big( (n-d) :(v^y)^{n-d-1} (v^h)^{d+1} C^{\beta bb}: -2 d : (v^y)^{n-d} (v^h)^{d-1} v^x C^{\beta bb}:+  :(v^y)^{n-d} (v^h)^d C^{\gamma bb}:  \big),$$ modulo terms which either have polynomial degree at most $2n+1$, or do not depend on $C^{\gamma bb}$. (Also, it is clear that the term $:(v^y)^{n-d-1} (v^h)^{d+1} C^{\beta bb}:$ cannot appear elsewhere in this expression).

Suppose first that $d<n$. The coefficient of $:(v^y)^{n-d} (v^h)^d C^{\gamma bb}:$ in (\ref{intercrap}) is $n-d+1+\lambda$, which is clearly nonzero, since $n-d\geq 1$ and $\lambda >-1$. Moreover, the coefficient of $:(v^y)^{n-d-1} (v^h)^{d+1} C^{\beta bb}:$ is $\lambda (n-d)$, and the ratio of these coefficients $\tilde{\lambda} =\frac{\lambda(n-d)}{n-d+1+\lambda}$ clearly lies in the interval $(-1,0)$, since $\lambda\in (-1,0)$. So we can divide (\ref{intercrap}) by $n-d+1+\lambda$, obtaining \begin{equation}\label{intercrapi}:(v^y)^{n-d} (v^h)^{d} C^{\gamma b b }:  + \tilde{\lambda} :(v^y)^{n-d-1} (v^h)^{d+1} C^{\beta bb}: -  \frac{2(d-1)}{n-d+1 +\lambda} :(v^y)^{n-d+1} (v^h)^{d-2}v^x C^{\gamma bb} : ,\end{equation} modulo terms which either have polynomial degree at most $2n+1$, or do not depend on $C^{\gamma bb}$. 

The expression (\ref{intercrapi}) is not quite of the desired form, because the term \begin{equation} \label{badterm}- \frac{2(d-1)}{n-d+1 +\lambda} :(v^y)^{n-d+1} (v^h)^{d-2}v^x C^{\gamma bb} : \end{equation} has polynomial degree $2n+3$ and depends on $C^{\gamma bb}$. However, by inductive assumption, there exists a polynomial of the form \begin{equation} \label{eliminator} :(v^y)^{n-d+1} (v^h)^{d-2}C^{\gamma bb}: + :f C^{\beta bb}: + \cdots \end{equation} in the kernel of $K(0)$, where $f$ has polynomial degree at most $2n-2$, and $(\cdots)$ has polynomial degree at most $2n-1$. Taking the Wick product of (\ref{eliminator}) with $\frac{2(d-1)}{n-d+1 +\lambda} v^x$, and adding it to (\ref{intercrapi}) eliminates the term (\ref{badterm}), but does not affect the coefficients of either $:(v^y)^{n-d} (v^h)^{d} C^{\gamma bb}:$ or $:(v^y)^{n-d-1} (v^h)^{d+1} C^{\beta bb}:$ in (\ref{intercrapi}). This yields a vertex operator of the form (\ref{specialvo}) in the case $d<n$.

Finally, suppose that $d=n$, which is easier than the case $d<n$ because there is no coefficient $\tilde{\lambda}$ to worry about. Then (\ref{intercrapi}) becomes
$$: (v^h)^{n} C^{\gamma b b }:  -  \frac{2(n-1)}{1 +\lambda} :(v^y) (v^h)^{n-2}v^x C^{\gamma bb} :,$$ modulo terms which either have polynomial degree at most $2n+1$, or do not depend on $C^{\gamma bb}$. As above, we can eliminate the term $-  \frac{2(n-1)}{1 +\lambda} :(v^y) (v^h)^{n-2}v^x C^{\gamma bb} :$ to obtain a vertex operator of the form (\ref{specialvo}). This completes the induction. \end{proof}

Consider the standard monomial basis for the universal enveloping algebra $U(\gs\gl_2)$ given by $$\{ x^r y^s h^t\in U(\gs\gl_2)| r,s,t\geq 0\}.$$ For any such monomial $\mu = x^r y^s h^t\in U(\gs\gl_2)$, we can associated to it the normally ordered monomial $:(v^x)^r(v^y)^s (v^h)^t: \in \cS^{\gs\gl_2[t]}$, which we also denote by $\mu$. Define $h_{\mu}\in \cW_{bas}\cap Ker(K(0))$ to be the Wick product of $:(v^x)^r:$ with the vertex operator $$: (v^y)^s (v^h)^t C^{\gamma bb} + :f C^{\beta bb}: + \cdots,$$ given by Lemma \ref{monomials}. Clearly $h_{\mu}$ has degree $4r-4s$, weight $2s+t+2$, and is of the form $$ :\mu C^{\gamma bb} : + :g C^{\beta bb}: + \cdots ,$$ where $g\in \cS^{\gs\gl_2[t]}$ has polynomial degree at most $2(r+s+t)$, and $(\cdots)$ has polynomial degree at most $2(r+s+t)+1$. In particular, for $\mu = 1$, we have $h_{\mu} = C^{\gamma bb}$. The assignment $\mu\mapsto h_{\mu}$ extends to a linear map $\phi: U(\gs\gl_2)\rightarrow \cW_{bas}^{(2)} \cap Ker(K(0))$, which induces a linear map $\Phi:U(\gs\gl_2)\rightarrow H^*(\cW_{bas},K(0))$. It is clear that $\Phi$ maps the eigenspace of $[h,-]$ of eigenvalue $d$ into $H^{2d}(\cW_{bas},K(0))$.

\begin{lemma} \label{lininjmap}The map $\Phi: U(\gs\gl_2)\rightarrow H^*(\cW_{bas},K(0))$ is injective. \end{lemma}

\begin{proof} Let $f$ be a nonzero, homogeneous element of degree $d$ in $U(\gs\gl_2)$. Then the leading term of $\phi(f)$ has polynomial degree $2d+3$, and is of the form $$:f C^{\gamma bb}: + :g C^{\beta bb}:,$$ for some $g\in \cS^{\gs\gl_2[t]}$ of polynomial degree at most $2d$. Since there are no nontrivial relations of the form (\ref{forbiddeni}), $\phi(f)$ is nonzero. Finally, since $\phi(f)$ has odd polynomial degree, it must represent a nontrivial class in $H^*(\cW_{bas},K(0))$. \end{proof}

\section{The structure of $\H^*_{SU(2)}(\mathbb{C})$}

In this section, we will use the spectral sequence of the double complex together with our results about $H^*(\cW_{bas},K(0))$ to study $\H^*_{SU(2)}(\mathbb{C})$. In particular, we will show that each nontrivial element of $H^*(\cW_{bas},K(0))$ given by Lemma \ref{lininjmap} gives rise to a nontrivial element of $\H^*_{SU(2)}(\mathbb{C})$. Recall that for each $k$, we have a decreasing filtration $$\H^{2k}_{SU(2)}(\mathbb{C})_{(0)} \supset  \H^{2k}_{SU(2)}(\mathbb{C})_{(2)} \supset   \H^{2k}_{SU(2)}(\mathbb{C})_{(4)} \supset \cdots,$$ where $\H^{2k}_{SU(2)}(\mathbb{C})_{(2r)}$ consists of classes admitting a representative of the form $\omega = \sum_{i\geq r} \omega^{(2i)}$, with $\omega^{(2i)}\in \cW_{bas}^{(2i)}$. We have an injective linear map \begin{equation}\label{lininjc}\H^{2k}_{SU(2)}(\mathbb{C})_{(2l)}/ \H^{2k}_{SU(2)}(\mathbb{C})_{(2l+2)} \rightarrow H^{2k}(\cW_{bas},J(0)),\end{equation} which is the specialization of (\ref{injbj}) to the case $G=SU(2)$.

\begin{lemma} \label{bnzero} Suppose that $\omega= \sum_{i\geq 0} \omega^{(2i)}$ is homogeneous with respect to degree and weight, with $\omega^{(2i)}\in \cW_{bas}^{(2i)}$, and $D(0)(\omega) = 0$. If $\omega^{(0)}\neq 0$, $\omega$ represents a nontrivial class in $\H^*_{SU(2)}(\mathbb{C})$. \end{lemma}

\begin{proof} This is clear from the injectivity of the maps (\ref{classinj}) and (\ref{lininjc}). \end{proof}

Using this result, there are two subalgebras of $\H^*_{SU(2)}(\mathbb{C})$ that we can now describe. The first one is essentially classical. Recall the abelian subalgebra $\langle \gamma \rangle^{\gs\gl_2[t]}$ of $\cS^{\gs\gl_2[t]}$, which is just the polynomial algebra with generators $\partial^k v^x$, $k\geq 0$. Since $D(0)$ acts trivially on $\langle \gamma\rangle ^{\gs\gl_2[t]}$, there is a vertex algebra homomorphism $\langle \gamma\rangle ^{\gs\gl_2[t]} \rightarrow \H^*_{SU(2)}(\mathbb{C})$, which is clearly injective by Lemma \ref{bnzero} and the fact that $\langle \gamma\rangle ^{\gs\gl_2[t]} \subset \cW_{bas}^{(0)}$. The next theorem is our main result; it describes an interesting and essentially nonclassical subalgebra of $\H^*_{SU(2)}(\mathbb{C})$.

\begin{thm} \label{main} There exists a linear map $\psi: U(\gs\gl_2)\rightarrow \cS^{\gs\gl_2[t]}$ such that for any $f\in U(\gs\gl_2)$, $\phi(f) + \psi(f)$ lies in $Ker(D(0))$. Moreover, the linear map $$\Psi: U(\gs\gl_2)\rightarrow \H^*_{SU(2)}(\mathbb{C})$$ sending $f$ to the class of $\phi(f) + \psi(f)$, is injective. \end{thm}

\begin{proof} Fix a nonzero monomial $\mu\in U(\gs\gl_2)$. Since $\phi(\mu) \in \cW_{bas}^{(2)}\cap Ker(K(0))$, and $J(0)$ preserves $\cW_{bas}$ and lowers $b$-number by $1$, we have $J(0)(\phi(\mu))\in \cW_{bas}^{(1)}$. Moreover, $J(0)(\phi(\mu)) \in Ker(K(0))$, since $K(0)$ commutes with $J(0)$. Since $\cW_{bas}^{(1)}\subset \cD$, and $H^{(1)}(\cD,K(0)) = 0$ by Lemma \ref{dcoh}, there an element $\omega_{\mu} \in \cD^{(0)} = \cS^{\gs\gl_2[t]}$ such that $K(0)(\omega_{\mu}) = -J(0)(\phi(\mu))$. Define the linear map $\psi: U(\gs\gl_2)\rightarrow \cS^{\gs\gl_2[t]}$ on our monomial basis by $\psi(\mu) = \omega_{\mu}$. Since $J(0)(\psi(\mu)) = 0$ by Lemma \ref{bnumbzero}, we clearly have $D(0)(\phi(\mu)+\psi(\mu)) = 0$. 

In order to prove that $\Psi$ is injective, it suffices to prove that $\psi$ is injective, by Lemma \ref{bnzero}. Let $f$ be a nonzero element of $U(\gs\gl_2)$ which is homogeneous of degree $d$. Then $\phi(f)$ has polynomial degree $2d+3$, and its leading term is of the form $$ : f C^{\gamma bb} : +: g C^{\beta bb}:$$ for some normally ordered polynomial $g$ in $v^x,v^y,v^h$ of polynomial degree at most $2d$. By Lemma \ref{bnumbzero}, $J(0)$ annihilates both $f$ and $g$, so by (\ref{opejci}) and (\ref{opejcii}) we have $$J(0)(\phi(f)) = 2: f v^h K: -4 : f v^x Q^{\beta b}: + 2: g v^y K:~+~ : g v^h Q^{\beta b}:,$$ modulo terms of polynomial degree at most $2d+1$. Since there are no relations of the form (\ref{forbiddenii}), the condition $J(0)(\phi(f)) = 0$ implies that, up to lower polynomial degree corrections, $$2 :f v^h: +2 :g v^y: = 0,\ \ \ \ \ \ \ \ -4 :f v^x: + :g v^h: = 0.$$ This in turn implies that $:\big(:v^h v^h:-4 :v^x v^y:\big) f: = 0$, up to lower polynomial degree corrections. This is impossible by Remark \ref{uca}, so $J(0)(\phi(f))\neq 0$. Since $K(0)(\psi(f)) = -J(0)(\phi(f))$, it follows that $\psi(f)\neq 0$ as well. \end{proof}

\begin{remark} It is clear from the construction that for $\mu = x^r y^s h^t$, $\Psi(\mu)$ has degree $4r-4s$ and conformal weight $2s+t+2$. For the identity element $1\in U(\gs\gl_2)$, we have \begin{equation} \label{defofl} \phi(1) + \psi(1)=C^{\gamma bb} + 2 :v^x v^y: + \frac{1}{2} :v^h v^h: -\frac{1}{2} \partial v^h.\end{equation} This coincides with the element $L$ given by (\ref{defofvirasoro}), expressed in terms of the generators of $\cC$. Hence $\Psi(1)$ is precisely the Virasoro class $\L\in \H^*_{SU(2)}(\mathbb{C})$. \end{remark}

For the sake of illustration, we write down explicit representatives for a few more of the classes in $\H^*_{SU(2)}(\mathbb{C})$ given by Theorem \ref{main}. For $n\geq 1$, let $H_{2n+2}$ denote $\phi(f) + \psi(f)$ for $f = y^n \in U(\gs\gl_2)$, which represents a class of degree $-4n$ and weight $2n+2$. Similarly, let $F_{n+2}$ denote $\phi(f)+\psi(f)$ for $f=h^n$, which represents a class of degree $0$ and weight $n+2$. The following formulae were verified using Kris Thielemann's Mathematica OPE package \cite{T}.
 
\begin{equation} F_3= :v^h C^{\gamma bb}: + \frac{2}{3} :v^x C^{\beta bb}:-\frac{5}{3} \partial C^{\gamma bb}\end{equation} $$+ \frac{4}{3} :v^y v^x v^h: + \frac{1}{3} :v^h v^h v^h: - \frac{1}{3}: v^h \partial v^h:  - 
 \frac{16}{3} :(\partial v^y)v^x: + \frac{2}{3} :v^y\partial v^x: - \frac{5}{3} \partial^2 v^h,$$

\begin{equation}F_4 = :v^h v^h C^{\gamma bb}:  +:( \partial v^h) C^{\gamma bb} :  +  : v^h v^x C^{\beta bb}: + 
 \frac{2}{3} :(\partial v^x) C^{\beta bb}: + \frac{1}{3}: v^x \partial C^{\beta bb}:\end{equation}  $$+  \frac{1}{4} :v^h v^h v^h v^h: +  :v^h v^h v^x v^y: + 
  \frac{1}{2}:v^h v^h \partial v^h:  + \frac{4}{3} :(\partial v^x) v^y v^h: + \frac{2}{3} v^x (\partial v^y) v^h: $$ $$ + \frac{5}{3} :v^x v^y \partial v^h: + 4 :(\partial^2 v^x)v^y:   - 2 :v^x \partial^2 v^y:  - \frac{1}{4} :(\partial v^h)(\partial v^h):  - \frac{1}{12} \partial^3 v^h,$$

\begin{equation}H_4 = ~:v^y C^{\gamma bb}:~ - \frac{1}{4} : v^h C^{\beta bb}:~  - \frac{5}{12} \partial C^{\beta bb}\end{equation} 
$$+  :v^x v^y v^y: + \frac{1}{4} :v^h v^h v^y: + \frac{7}{6} :v^h \partial v^y: - 
 \frac{19}{12} :(\partial v^h) v^y: + \frac{1}{12} \partial^2 v^y,$$
 
\begin{equation}H_6 = :v^y v^y C^{\gamma bb}:  - \frac{1}{3} :v^y v^h C^{\beta bb}: -  \frac{2}{15} :(\partial v^y) C^{\beta bb}:  - \frac{7}{15} :v^y \partial C^{\beta bb}:\end{equation}  $$ + \frac{1}{6} :v^y v^y v^h v^h: + \frac{2}{3}: v^x v^y v^y v^y: + \frac{2}{15}:(\partial v^y) v^y v^h: - \frac{53}{30} :v^y v^y \partial v^h: + 2 :v^y \partial^2 v^y:,$$
 
\begin{equation}H_8 = :v^y v^y v^y C^{\gamma bb}: - \frac{3}{8} :v^y v^y v^h C^{\beta bb}: - 
 \frac{3}{14} :(\partial v^y) v^y C^{\beta bb}:  -  \frac{27}{56} :v^y v^y \partial C^{\beta bb}: \end{equation}  $$+ 
\frac{1}{2} : v^x v^y v^y v^y v^y: + \frac{1}{8}:v^y v^y v^y v^h v^h: - 
 \frac{103}{56} :v^y v^y v^y \partial v^h: + \frac{3}{28} :(\partial v^y) v^y v^y v^h:  + 
 3 :(\partial^2 v^y) v^y v^y:. $$
 
We conjecture that the classes $\{\Psi(f)|~ f\in U(\gs\gl_2)\}$, together with the classical element $[v^x]\in \H^*_{SU(2)}(\mathbb{C})[0] \cong H^*_{SU(2)}(pt)$ represented by $v^x$, form a strong generating set for $\H^*_{SU(2)}(\mathbb{C})$. If this is the case, $\H^*_{SU(2)}(\mathbb{C})$ has purely even degree. In fact, it is likely that $\H^*_{SU(2)}(\mathbb{C})$ is generated (although not strongly generated) by a much more economical set. The following circle product relations have been verified:
\begin{equation} \label{somecircprod}H_4\circ_1 H_4 = \frac{5}{2}H_6 , \ \ \ \ \ \ \ H_4\circ_1 H_6 = \frac{112}{45} H_8, \ \ \ \ \ \ \  H_4\circ_1 v^x = -\frac{5}{12} L.\end{equation} We expect that there exist nonzero constants $c_n$ such that $H_4\circ_1 H_{2n+2} = c_n H_{2n+4}$. If this holds, all the $H_{2n+2}$ lie in the vertex algebra generated by $H_4$. It is clear from the proof of Lemma \ref{monomials} that all elements of the form $\phi(f)$ for $f\in U(\gs\gl_2)$ lie in the vertex algebra generated by $v^x$ and $h_{2n+2}$ for $n\geq 2$. It follows that the classes $\Psi(f)\in \H^*_{SU(2)}(\mathbb{C})$ for $f\in U(\gs\gl_2)$ all lie in the vertex algebra generated by $[v^x]$ and $[H_{2n+2}]$ for $n\geq 2$. (Here $[H_{2n+2}]$ denotes the class represented by $H_{2n+2}$). If all the $H_{2n+2}$ lie in the vertex algebra generated by $H_4$ as suggested by (\ref{somecircprod}), each $\Psi(f)$ lies in the vertex algebra generated by $[v^x]$ and $[H_4]$. Finally, if $\H^*_{SU(2)}(\mathbb{C})$ is indeed strongly generated by $\{[v^x], \Psi(f)|~f\in U(\gs\gl_2)\}$, this would imply that $\H^*_{SU(2)}(\mathbb{C})$ is generated as a vertex algebra by $[v^x]$ and $[H_4]$.

An interesting problem is to compute the graded character 
$$\chi(z,q) = \sum_{d\in\mathbb{Z}} \sum_{n\geq 0} \text{dim}~ \H^d_{SU(2)}(\mathbb{C})[n] z^d q^n,$$ where  $\H^d_{SU(2)}(\mathbb{C})[n]$ has degree $d$ and weight $n$. Even if $\H^*_{SU(2)}(\mathbb{C})$ is strongly generated by $\{[v^x], \Psi(f)|~f\in U(\gs\gl_2)\}$, it is not clear how to compute this character, because there exist normally ordered polynomial relations among these generators and their derivatives. For example, in weight $4$ and degree $0$, a computer calculation shows that the following expression in $\cW_{bas}$ is identically zero: $$:L L: -F_4  - 4 :v^x H_4: +  \partial F_3+ \frac{7}{6} \partial^2 L.$$ 

Finally, we expect that there exists an alternative and more geometric construction of the chiral equivariant cohomology, and in particular of $\H^*_G(\mathbb{C})$. We hope that the structure described in this paper in the case $G=SU(2)$ may give some hint about where to look for such a construction. A natural place to begin would be to find a geometric interpretation of the element $[H_4]\in \H^{-4}_{SU(2)}(\mathbb{C})[4]$.

\section{Appendix}
In this Appendix we develop some techniques for studying invariant rings of the form $\cO(J_{\infty}(V))^{\gg[t]}$ and we prove Theorem \ref{jet}. \footnote[1]{The localization technique used to study invariant rings of the form $\cO(J_m(V))^{\gg[t]}$ in this section is due to Bailin Song. I thank him for sharing this idea with me.} First, we recall some basic facts about jet schemes, following the notation in \cite{EM}. Let $X$ be an irreducible scheme of finite type over $\mathbb{C}$. For each integer $m\geq 0$, the jet scheme $J_m(X)$ is determined by its functor of points: for every $\mathbb{C}$-algebra $A$, we have a bijection
$$Hom (Spec (A), J_m(X)) \cong Hom (Spec (A[t]/\langle t^{m+1}\rangle ), X).$$ Thus the $\mathbb{C}$-valued points of $J_m(X)$ correspond to the $\mathbb{C}[t]/\langle t^{m+1}\rangle$-valued points of $X$. If $m>p$, we have projections $\pi_{m,p}: J_m(X) \rightarrow J_p(X)$ which are compatible when defined: $\pi_{m,p} \circ \pi_{q,m} = \pi_{q,p}$. Clearly $J_0(X) = X$ and $J_1(X)$ is the total tangent space $Spec(Sym(\Omega_{X/\mathbb{C}}))$. The assignment $X\mapsto J_m(X)$ is functorial, and a morphism $f:X\ra Y$ of schemes induces $f_m: J_m(X) \ra J_m(Y)$ for all $m\geq 1$. If $X$ is nonsingular, $J_m(X)$ is irreducible and nonsingular for all $m$. Moreover, if $X,Y$ are nonsingular and $f:Y\ra X$ is a smooth surjection, $f_m$ is surjective for all $m$. 

If $X=Spec(R)$ where $R= \mathbb{C}[y_1,\dots,y_r] / \langle f_1,\dots, f_k\rangle$, we can find explicit equations for $J_m(X)$. Define new variables $y_j^{(i)}$ for $i=0,\dots, m$, and define a derivation $D$ on the generators of $\mathbb{C}[y_1^{(i)},\dots, y_r^{(i)}]$ by $D(y_j^{(i)}) = y_j^{(i+1)}$ for $i<m$, and $D(y_j^{(m)}) =0$. Since $D$ is a derivation, this specifies the action of $D$ on all of $\mathbb{C}[y_1^{(i)},\dots, y_r^{(i)}]$; in particular, $f_j^{(i)} = D^i ( f_j)$ is a well-defined polynomial in $\mathbb{C}[y_1^{(i)},\dots, y_r^{(i)}]$. Letting $R_m = \mathbb{C}[y_1^{(i)},\dots, y_r^{(i)}] / \langle f_1^{(i)},\dots, f_k^{(i)}\rangle$, we have $J_m(X)\cong Spec (R_m)$. 

Given a scheme $X$, define $J_{\infty}(X) = \lim_{\infty \leftarrow m} J_m(X)$, which is known as the infinite jet scheme, or space of arcs of $X$. If $X = Spec(R)$ as above, $J_{\infty}(X)\cong Spec(R_{\infty})$ where $R_{\infty}  = \mathbb{C}[y_1^{(i)},\dots, y_r^{(i)}] / \bra f_1^{(i)},\dots, f_k^{(i)}\ket$. Here $i=0,1,2,\dots$ and $D (y^{(i)}_j) = y^{(i+1)}_j$ for all $i$. We denote by $\cO(J_{\infty}(X))$ the ring $\lim_{m\rightarrow \infty} \cO(J_m(X))$.

Let $G$ be a connected, reductive algebraic group over $\mathbb{C}$ with Lie algebra $\gg$. For $m\geq 1$, $J_m(G)$ is an algebraic group which is the semidirect product of $G$ with a unipotent group $U_m$. The Lie algebra of $J_m(G)$ is $\gg[t]/t^{m+1}$. Given a linear representation $V$ of $G$, there is an action of $G$ on $\cO(V)$ by automorphisms, and a compatible action of $\gg$ on $\cO(V)$ by derivations, satisfying $$\frac{d}{dt} exp (t\xi) (f)|_{t=0} = \xi(f),\ \ \ \ \xi\in\gg,\ \ \ \ f\in \cO(V).$$ Choose a basis $\{x_1,\dots,x_n\}$ for $V^*$, so that $$\cO(V) \cong  \mathbb{C}[x_1,\dots,x_n],\ \ \ \ \ \ \cO(J_m(V)) =  \mathbb{C}[x_1^{(i)},\dots,x_n^{(i)}|~,0\leq i\leq m].$$ We regard $\cO(V)$ as a subalgebra of $\cO(J_m(V))$ by identifying $x_i$ with $x_i^{(0)}$. Then $J_m(G)$ acts on $J_m(V)$, and the induced action of $\gg[t]/t^{m+1}$ by derivations on $\cO(J_m(V))$ is defined on generators by \begin{equation}\label{jetaction}\xi t^r (x_j^{(i)}) = c^r_i \xi(x_j)^{(i-r)},\end{equation} where $c^r_i = \frac{i!}{(i-r)!}$ for $0\leq r\leq i$, and $c^r_i = 0$ for $r>i$. Via the projection $\gg[t]\ra \gg[t]/t^{m+1}$, $\gg[t]$ acts on $\cO(J_m(V))$, and the invariant rings $\cO(J_m(V))^{\gg[t]}$ and $\cO(J_m(V))^{\gg[t]/t^{m+1}}$ coincide.

The problem of describing invariant rings of the form $\cO(J_m(V))^{\gg[t]}$ was first studied (to the best of our knowledge) by D. Eck in \cite{E}. Eck was primarily interested in the case where $G$ and $V$ are real, although in \cite{E} he worked with the complexifications of $G$ and $V$ in order to use tools from algebraic geometry. He proved that for all $m\geq 1$, \begin{equation} \label{eckequation} \cO(J_m(V))^{\gg[t]} = \langle \cO(V)^G\rangle_m\end{equation} under fairly restrictive hypotheses, namely, that the categorical quotient $V//G$ is smooth and the map $V\rightarrow V//G$ is an \'etale fibration. In this notation, $\langle \cO(V)^G\rangle_m$ denotes the ring generated by $\{D^i (f)|~f\in \cO(V)^G|~0\leq i\leq m\}$. If (\ref{eckequation}) holds, letting $m$ approach infinity shows that $\cO(J_{\infty} (V))^{\gg[t]} = \langle \cO(V)^G\rangle$, which is just the ring generated by $\{D^i (f) |~f\in \cO(V)^G,~ i\geq 0\}$. In the case where $V$ is the adjoint representation of a simple group $G$, this result was also proven in the appendix of \cite{Mu}. More recently \cite{P}, the equality $\cO(J_{\infty} (V))^{\gg[t]} = \langle \cO(V)^G\rangle$ was proven in the case where $G$ is a compact, real Lie group and $V$ is a real, irreducible representation of $G$. 

For the moment we make no additional assumptions about $G$ and $V$. Since $G$ is reductive, $\cO(V)^G$ is finitely generated. Choose a set of irreducible, homogeneous generators $\{y_1,\dots, y_p\}$ for $\cO(V)^G$. Let $d = \text{dim}(V//G)$, and let \begin{equation}\label{defofvo}V^{0}= \{x\in V|~\text{rank} \big[\frac{\partial y_i}{\partial x_j}\big]\big|_{x} = d\},\end{equation} which is $G$-invariant, Zariski open, and independent of our choice of generators for $\cO(V)^G$. For each $x\in V^{0}$, we can choose $d$ algebraically independent polynomials from the set $\{y_1,\dots,y_p\}$ such that the corresponding $d\times n$ submatrix of $\big[\frac{\partial y_i}{\partial x_j}\big]$ has rank $d$. Without loss of generality, we may assume these are $y_1,\dots, y_d$. It is easy to see that the polynomials $\{D^i (y_1),\dots, D^i (y_d)|~ i\geq 0\}$ are algebraically independent, and distinct monomials in the variables $D^i ( y_j)$ are linearly independent over the function field $\cK(V)= \mathbb{C}(x_1,\dots, x_n)$.

For notational simplicity, we will denote $\cO(V)$ by $\cO$, and we will denote $\cO(J_m(V))$ by $\cO_m$. The ring $\cO_m$ has a $\mathbb{Z}_{\geq 0}$-grading $ \cO_m = \bigoplus_{r\geq 0} \cO_m[r]$ by weight, defined by $wt(x_i^{(j)}) = j$. Note that $\cO_m[0] = \cO$ for all $m$, and each $\cO_m[r]$ is a module over $\cO$. By (\ref{jetaction}), for each $\xi\in\gg$, the action of $\xi t^k$ is homogeneous of weight $-k$. In particular, the Lie subalgebra $t\gg[t]\subset \gg[t]$ annihilates $\cO$, and hence acts by $\cO$-linear derivations on each $\cO_m$.

Let $\cK$ be the quotient field of $\cO$, and let $$\cK_m  = \cK \otimes_{\cO} \cO_m,\ \ \ \ \ \cK_m[r] = \cK\otimes_{\cO} \cO_m[r].$$ Clearly $t\gg[t]$ acts on $\cK_m$ by $\cK$-linear derivations. Our first goal is to describe the invariant space $\cK_1^{t\gg[t]}$.

\begin{lemma} \label{ckinvariants} $\cK_1^{t\gg[t]}$ is generated by $y_1^{(1)},\dots,y_d^{(1)}$ as a $\cK$-algebra.\end{lemma}

\begin{proof} Let $W\subset \cK_1[1]$ be the $\cK$-subspace spanned by $y_1^{(1)},\dots, y_d^{(1)}$, which has dimension $d$ over $\cK$. Since $G$ is reductive and acts on $\cK_1[1]$, $W$ has a $G$-stable complement $A\subset \cK_1[1]$ of dimension $r=n-d$. The linear map $$\gg\ra Hom(A,\cK),\ \ \ \ \ \ \xi\mapsto \xi t|_A$$ is surjective, so we can choose $\xi_1,\dots,\xi_r\in \gg$ such that $\{\xi_1 t|_A,\dots, \xi_r t|_A\}$ form a basis for $Hom(A,\cK)$. Choose a dual basis $\{a_1,\dots, a_r\}$ for $A$ satisfying \begin{equation} \label{duality}\xi_i t (a_j) = \delta_{i,j}.\end{equation}
Since $\{y_1^{(1)},\dots,y_d^{(1)}, a_1,\dots, a_r\}$ is a basis for $\cK_1[1]$ over $\cK$, it follows that $$\cK_1 = \cK \otimes_{\cO} Sym (\cK_1 [1]) = \cK \otimes_{\cO} Sym(W)\otimes_{\cO} Sym(A).$$
Since $y_i^{(1)}\in \cK_1[1]^{t\gg[t]}$ and each $\xi t$ is a $\cK$-linear derivation on $\cK_1$, it follows that $Sym(W)^{t\gg[t]} = Sym(W)$. Similarly, (\ref{duality}) shows that $Sym(A)^{t\gg[t]} = \mathbb{C}$. Hence $\cK_1^{t\gg[t]} = \cK \otimes_{\cO} Sym(W)$, as claimed. \end{proof}

In fact, it suffices to work over a certain localization of $\cO$ rather than the full quotient field $\cK$. We can choose $r=n-d$ elements of the set $\{x_1^{(1)},\dots, x_n^{(1)}\}$, say $x_{d+1}^{(1)},\dots, x_n^{(1)}$, such that 
\begin{equation} \label{basiski} \{y_1^{(1)},\dots, y_d^{(1)}, x_{d+1}^{(1)},\dots, x_n^{(1)}\}\end{equation} forms a basis for $\cK_1[1]$ over $\cK$.  It follows that the set (\ref{basiski}) is algebraically independent. Let $\Delta$ be the determinant of the $\cK$-linear change of coordinates on $\cK_1[1]$ given by \begin{equation}\label{flinear}(x_1^{(1)},\dots,x_d^{(1)}, x_{d+1}^{(1)},\dots,x_n^{(1)})\mapsto (y_1^{(1)},\dots, y_d^{(1)}, x_{d+1}^{(1)},\dots,x_n^{(1)}),\end{equation} and let $\cO_{\Delta}$ be the localization of $\cO$ along the multiplicative set generated by $\Delta$. Let $$\cO_{m,\Delta} = \cO_{\Delta} \otimes_{\cO} \cO_m,$$ which has a weight grading $ \cO_{m,\Delta} = \bigoplus_{k\geq 0} \cO_{m,\Delta}[k]$.

\begin{lemma} \label{levelonecase}The invariant space $(\cO_{1,\Delta})^{t\gg[t]}$ is generated as an $\cO_{\Delta}$-algebra by $y_1^{(1)},\dots,y_d^{(1)}$.\end{lemma}

\begin{proof} Any $\omega\in (\cO_{1,\Delta})^{t\gg[t]}$ can be expressed uniquely in the form \begin{equation}\label{monomialsumi}\omega = \sum_{i\in I} p_i\mu_i,\end{equation} where $\mu_i$ are distinct monomials in $y_1^{(1)},\dots,y_d^{(1)}$ and $p_i$ are polynomials in $x_{d+1}^{(1)},\dots,x_n^{(1)}$ with coefficients in $\cO_{\Delta}$.

Since $t\gg[t]$ acts trivially on each $\mu_i$, and $\{y_1^{(1)},\dots,y_d^{(1)}, x_{d+1}^{(1)},\dots,x_n^{(1)}\}$ are algebraically independent, it follows that each $p_i$ lies in $(\cO_{1,\Delta})^{t\gg[t]}$ independently. Therefore we may assume without loss of generality that $\omega=p(x_{d+1}^{(1)},\dots,x_n^{(1)})$. But by Lemma \ref{ckinvariants}, $\omega$ can also be expressed as a polynomial $p'(y_1^{(1)},\dots,y_d^{(1)})$ with $\cK$-coefficients. Thus $$p(x_{d+1}^{(1)},\dots,x_n^{(1)}) = p'(y_1^{(1)},\dots,y_d^{(1)}),$$ which violates the algebraic independence of $\{ y_1^{(1)},\dots, y_d^{(1)}, x_{d+1}^{(1)},\dots,x_n^{(1)}\}$. \end{proof}

\begin{lemma} \label{tgtinvariants} For each $m>0$, $(\cO_{m,\Delta})^{tg[t]}$ is generated as an $\cO_{\Delta}$-algebra by $$\{y_1^{(j)},\dots, y_d^{(j)}|~1\leq j\leq m\}.$$\end{lemma} 
\begin{proof} 
This holds for $m=1$ by Lemma \ref{levelonecase}, so we may assume inductively that $(\cO_{m-1,\Delta})^{tg[t]}$ is generated as an $\cO_{\Delta}$-algebra by $\{y_1^{(j)},\dots, y_d^{(j)}|,~1\leq j<m\}$. Let $I\subset \cO_{m,\Delta}$ denote the ideal generated by $\{x_1^{(m)},\dots, x_n^{(m)}\}$, and consider the filtration $$ \cO_{m,\Delta} \supset I \supset I^2 \supset\cdots ,$$ and the associated grading $\cO_{m,\Delta} \cong  \bigoplus_{k\geq 0} I^k/ I^{k+1}$. Let $S=\mathbb{C}[x_i, x_i^{(m)}]$ and let $S_{\Delta} = S\otimes_{\cO} \cO_{\Delta}$. We have decompositions \begin{equation}\label{decompom}\cO_{m} = S\otimes_{\cO} \cO_{m-1},\ \ \ \ \ \ \cO_{m,\Delta} = S_{\Delta} \otimes_{\cO} \cO_{m-1,\Delta}.\end{equation}  
Given $\omega\in (\cO_{m,\Delta})^{t\gg[t]}$, we say that $\omega$ has $m$-degree $s$ if $s$ is the minimal integer for which $\omega \in \bigoplus_{k=0}^s I^k/ I^{k+1}$. If $\omega$ has $m$-degree $s$, let $\hat{\omega}$ be the \lq\lq leading term" of $\omega$, i.e., the projection of $\omega$ onto $I^s/I^{s+1}$. Since $t^m\gg[t]$ acts trivially on $\cO_{m-1,\Delta}$, it is immediate from the decomposition (\ref{decompom}) that $$\hat{\omega} \in (S_{\Delta})^{t^m \gg[t]} \otimes_{\cO} ( \cO_{m-1,\Delta})^{t\gg[t]}.$$
There is a natural map of $\cO_{\Delta}$-algebras $S_{\Delta} \ra \cO_{1,\Delta}$, defined on generators by \begin{equation}\label{mtoone} x_j^{(m)} \mapsto x_j^{(1)},\ \ \ \ \ \ x_j \mapsto x_j.\end{equation} Moreover, given $\xi\in\gg$, $$\xi t^m (x_j^{(m)}) = m! \xi t (x_j^{(1)}) = m! \xi(x_j).$$  Hence $(S_{\Delta})^{t^m \gg[t]}\cong (\cO_{1,\Delta})^{t\gg[t]}$ as $\cO_{\Delta}$-algebras. 

For $i=1,\dots,d$, the leading term $\hat{y}_i^{(m)}$ of $y_i^{(m)}$ lies in $I/ I^2$, and $$y_i^{(m)} - \hat{y}_i^{(m)}\in I^0 = \cO_{m-1,\Delta}.$$ Since $\hat{y}_i^{(m)} \mapsto y_i^{(1)}$ under (\ref{mtoone}), it follows from Lemma \ref{levelonecase} that $(S_{\Delta})^{t^m \gg[t]}$ is generated by $\{\hat{y}_i^{(m)}|~i=1,\dots,d\}$ as an $\cO_{\Delta}$-algebra. 

By our inductive assumption, $(\cO_{m-1,\Delta})^{t\gg[t]}$ is generated by $\{y_i^{(j)}|~1\leq j <m\}$ as an $\cO_{\Delta}$-algebra. Since $\hat{\omega}\in (S_{\Delta})^{t^m \gg[t]} \otimes ( \cO_{m-1,\Delta})^{t\gg[t]}$, it follows that $\hat{\omega}$ can be expressed as a polynomial in the variables $$\{\hat{y}_i^{(m)},~y_i^{(j)}|,~1\leq j<m\}$$ with coefficients in $\cO_{\Delta}$. Moreover, $\hat{\omega}$ is homogeneous of degree $s$ in the variables $\hat{y}_i^{(m)}$ since $\hat{\omega}\in I^s/I^{s+1}$. Let $\omega'$ be the polynomial in the variables $\{y_i^{(j)}|~0\leq j\leq m\}$ obtained from $\hat{\omega}$ by replacing the variables $\hat{y}_i^{(m)}$ with $y_i^{(m)}$. Clearly $$\omega'' = \omega - \omega' \in (\cO_{m,\Delta})^{t\gg[t]},$$ and $\omega''$ has $m$-degree at most $s-1$. Since $\omega \sim \omega''$ modulo the $\cO_{\Delta}$-algebra generated by $\{y_i^{(j)}|~ 1\leq j\leq m\}$, the claim follows by induction on $m$-degree. \end{proof}

At this point, we impose a mild technical condition. We say that $\cO$ {\it contains no invariant lines} if $\cO$ contains no nontrivial, one-dimensional $G$-invariant subspace. If $G$ is semisimple, this condition is automatically satisfied by any $V$.

\begin{lemma} \label{charfree}If $V$ is a representation of $G$ for which $\cO$ contains no invariant lines, $\cO^G$ is generated by primes.\end{lemma}

\begin{proof} Let $f\in \cO^G$, and suppose that $f = p_1\cdots p_k$ is the prime factorization of $f$ in $\cO$. Since $\cO$ contains no invariant lines, each $g\in G$ must permute the factors of $f$. Hence $f$ determines a group homomorphism $\phi: G\ra S_k$ where $S_k$ is the permutation group on $k$ letters. But $G$ is connected and $\phi$ is continuous, so $\phi$ is trivial, and each $p_i\in\cO^G$. \end{proof}

Suppose that $\cO$ contains no invariant lines, and consider the full invariant algebra $\cO_m^{\gg[t]}$, which is just the $G$-invariant subalgebra $(\cO_m^{t\gg[t]})^G$. Given $\omega\in \cO_m^{\gg[t]}$, we may write $\omega = \sum_k p_k \mu_k$ where $\mu_k$ are distinct monomials in $\{y_1^{(j)},\dots,y_d^{(j)}|~ j=1,\dots,m\}$, and $p_k\in \cO_{\Delta}$. Since each $\mu_k\in \cO_m^{\gg[t]}$, and the $\mu_k$'s are linearly independent over $\cO_{\Delta}$, it follows that each $p_k\in \cO_{\Delta}^G$. If we express $p_k$ as a rational function $\frac{f}{g}$ in lowest terms, the denominator will either be 1, or will contain only $G$-invariant prime factors of $\Delta$, since both the denominator and numerator must be invariant. Letting $\Delta'$ be the product of the $G$-invariant prime factors of $\Delta$, it follows from Lemma \ref{tgtinvariants} that \begin{equation}\label{localization} \cO_{\Delta'}^{G} \otimes_{\cO^{G}} \cO_m^{\gg[t]} = \cO_{\Delta'}^{G} \otimes_{\cO^{G}} \bra \cO^{G}\ket_m,\end{equation} where $\bra \cO^{G}\ket_m\subset \cO_m^{\gg[t]}$ is generated by $\{D^{j} (f)|~ f\in\cO^G,~ 0\leq j\leq m\}$. If $\Delta$ contains no $G$-invariant prime factors (so that $\Delta' = 1$), it is immediate that $\cO_m^{\gg[t]}= \langle \cO^G \rangle_m$ for all $m\geq 1$. 

Even if $\Delta'\neq 1$, it still may be possible to prove the equality $\cO_m^{\gg[t]}=\bra \cO^{G}\ket_m$ for $m\geq 1$ using this method. Recall that $\Delta'$ corresponded to a choice of algebraically independent elements $\{y_1,\dots,y_d\}\subset \cO^G$. For $i=1,\dots,r$, let $\{y^i_1,\dots,y^i_d\}$ be a collection of maximal algebraically independent subsets of $\cO^G$, and let $\Delta'_1,\dots,\Delta'_r$ be the corresponding elements of $\cO^G$, obtained as above. It is easy to see from Lemma \ref{charfree} that in order to prove the equality $\cO_m^{\gg[t]}=\bra \cO^{G}\ket_m$ for all $m\geq 1$, it suffices to show that $gcd(\Delta'_1,\dots,\Delta'_r) = 1$.

{\it Proof of Theorem \ref{jet}.} First we consider the case of two copies of the adjoint representation $V$. As in Theorem \ref{weyl}, we work in the basis $\{a_i^x, a^y_i, a^h_i|~i=1,2\}$. Recall that $\cO(V\oplus V)^G$ has generators $q_{11}$, $q_{12}$, and $q_{22}$, which are algebraically independent. Consider the change of variables $$\big((a_1^h)^{(1)}, (a_2^h)^{(1)}, (a_1^x)^{(1)}, (a_2^x)^{(1)}, (a_1^y)^{(1)},(a_2^y)^{(1)}\big) \mapsto \big(q_{11}^{(1)},q_{12}^{(1)},q_{22}^{(1)}, (a_2^x)^{(1)},(a_1^y)^{(1)},(a_2^y)^{(1)}\big).$$ The determinant $\Delta$ of this change of variables is $8 a_2^h (a_1^y a_2^h - a_1^h a_2^y)$, so $\Delta' =1$. This proves (\ref{jettwo}).

Next, we consider the case of three copies of $V$, and we work in the basis $\{a_i^x, a^y_i, a^h_i|~i=1,2,3\}$. The ring $\cO(V\oplus V\oplus V)^{G}$ has generators $q_{11}$, $q_{12}$, $q_{13}$, $q_{22}$, $q_{23}$, $q_{33}$, and $c_{123}$, and the ideal of relations is generated by $$(c_{123})^2 +\frac{1}{4}\left| \begin{array}{lll} q_{11} & q_{12} &
q_{13} \\q_{12} & q_{22} & q_{23} \\ q_{13} & q_{23} & q_{33}\end{array}\right|.
$$ Clearly $(V\oplus V\oplus V)//G$ has dimension $6$. First, take $\{y_1,y_2,y_3,y_4,y_5,y_6\} = \{q_{11},q_{12},q_{22},q_{33},q_{23},q_{13}\}$, which is clearly algebraically independent, and consider the change of variables 
$$\big((a^h_1)^{(1)},(a^h_2)^{(1)},(a^h_3)^{(1)},(a^x_1)^{(1)},(a^x_2)^{(1)}, (a^y_1)^{(1)},(a^x_3)^{(1)}, (a^y_2)^{(1)}, (a^y_3)^{(1)}\big) $$ $$ \mapsto \big(q_{11}^{(1)},q_{12}^{(1)},q_{22}^{(1)},q_{33}^{(1)},q_{23}^{(1)},q_{13}^{(1)}, (a^x_3)^{(1)}, (a^y_2)^{(1)}, (a^y_3)^{(1)}\big).$$ We have $\Delta = 64 a^h_3 (-a_3^y a_2^h + a_3^h a_2^y) c_{123}$, so $\Delta' = c_{123}$.

Now consider $\{y_1,y_2,y_3,y_4,y_5,y_6\} = \{c_{123},q_{12},q_{22},q_{33},q_{23},q_{13}\}$, and consider the change of variables 
$$\big((a^h_1)^{(1)},(a^h_2)^{(1)},(a^h_3)^{(1)},(a^x_1)^{(1)},(a^x_2)^{(1)}, (a^y_1)^{(1)},(a^x_3)^{(1)}, (a^y_2)^{(1)}, (a^y_3)^{(1)}\big) $$ $$\mapsto (c_{123}^{(1)},q_{12}^{(1)},q_{22}^{(1)},q_{33}^{(1)},q_{23}^{(1)},q_{13}^{(1)}, (a^x_3)^{(1)}, (a^y_2)^{(1)}, (a^y_3)^{(1)}\big).$$ We have $\Delta = -8 a_3^h (-a_3^y a_2^h + a_3^h a_2^y)(q_{22} q_{33} - q_{23} q_{23}) $, so $\Delta' = q_{22} q_{33} - q_{23} q_{23}$. Since $c_{123}$ and $q_{22} q_{33} - q_{23} q_{23}$ have no common factor, (\ref{jetthree}) follows. $\Box$

\end{document}